\documentclass[twoside]{article}

\usepackage[accepted]{aistats2021}
\usepackage{lmodern}
\usepackage[utf8]{inputenc} 
\usepackage[T1]{fontenc}    
\usepackage{hyperref}       
\usepackage{url}            
\usepackage{booktabs}       
\usepackage{amsfonts}       
\usepackage{nicefrac}       
\usepackage{microtype}      
\usepackage{amsmath}
\usepackage{amsthm}
\usepackage{amssymb}
\usepackage{algorithmic}
\usepackage{graphicx}
\usepackage[ruled,vlined]{algorithm2e}
\usepackage{bbm}
\usepackage{mathrsfs}
\usepackage[round]{natbib}
\newtheorem{theorem}{Theorem}
\newtheorem{definition}{Definition}
\newtheorem{lemma}{Lemma}
\newtheorem{remark}{Remark}
\newtheorem{proposition}{Proposition}
\newtheorem{assumption}{Assumption}
\newtheorem{corollary}{Corollary}
\newcommand{\norm}[1]{\left\| #1 \right\|}
\usepackage{color}

\begin{document}

%

%

\twocolumn[
\aistatstitle{Reinforcement Learning for Mean Field Games with Strategic Complementarities}
\aistatsauthor{Kiyeob Lee\textsuperscript{*} \And Desik Rengarajan\textsuperscript{*}  \And Dileep Kalathil \And Srinivas Shakkottai}
\aistatsaddress{ Texas A\&M University  \\ \texttt{\{kiyeoblee,desik,dileep.kalathil,sshakkot\}@tamu.edu} } 
]

\begin{abstract}
\let\thefootnote\relax\footnotetext{\textsuperscript{*}Equal Contribution}
Mean Field Games (MFG) are the class of games with a very large number of agents and the standard equilibrium concept  is a Mean Field Equilibrium (MFE). Algorithms for learning MFE in dynamic MFGs are unknown in general.  Our focus is on an important subclass that possess a monotonicity property called Strategic Complementarities (MFG-SC).  We introduce a natural refinement to the equilibrium concept that we call Trembling-Hand-Perfect MFE (T-MFE),  which allows agents to employ a measure of randomization while accounting for the impact of such randomization on their payoffs.   We propose a simple algorithm for computing T-MFE under a known model.  We also introduce  a model-free and a model-based approach to learning T-MFE  and provide sample complexities of both algorithms.  We also develop a fully online learning scheme that obviates the need for a simulator.  Finally, we empirically evaluate the performance of the proposed algorithms via examples motivated by real-world applications.
\end{abstract}

\section{Introduction}

\emph{Strategic complementarities} refers to a well-established strategic game structure wherein the marginal returns to increasing one's strategy rise with increases in the competitors' strategies, and many qualitative results of the dynamic systems are based on these properties~\citep{milgrom1990rationalizability,vives2009strategic,adlakha2013mean}. Many practical scenarios with this property, including pricing in oligopolistic markets, adoption of network technology and standards, deposits and withdrawals in banking, weapons purchasing in arms races etc., have been identified, and our running example is that of adopting actions to prevent the spread of a computer or human virus, wherein stronger actions towards maintaining health (eg., installing patches, or wearing masks) by members enhances the returns (eg., system reliability or economic value) to a particular individual following suit.  

The above examples are characterized by a large number of agents following the dynamics of repeated action, reward, and state transition that is a characteristic of stochastic games.  However, analytical complexity implies that most work has focused on the static scenario with a small number of agents (see~\cite{milgrom1990rationalizability}).  Recently, there have been attempts to utilize the information structure of a mean field game (MFG)~\citep{lasry2007mean,tembine2009mean,adlakha2013mean,iyer2014mean,li2016incentivizing,LiXia18} to design algorithms to compute equilibria under a known model in the large population setting~\citep{adlakha2013mean}.  Here, each agent assumes that the states of all others are drawn in an i.i.d. manner from a common mean field belief distribution, and optimizes accordingly. However, the model is non-stationary due to the change in the mean field distribution at each time step, and provably convergent learning algorithms for identifying equilibrium strategies are currently unavailable. 

\paragraph{Main Contributions:} We study the problem of learning under an unknown model in stochastic games with strategic complementarities under the mean field setting.  Our main contributions are: \vspace{0.05in}\\
(i) We introduce the notion of trembling-hand perfection to the context of mean field games, under which a known randomness is introduced into all strategies.  Unlike an $\epsilon$-greedy policy in which randomization is added as an afterthought to the optimal action, under trembling-hand perfection, optimal value is computed consistently while accounting for this randomness.\vspace{0.05in}\\
(ii) We describe TQ-value iteration, based on generalized value iteration, as a means of computing trembling-hand consistent TQ-values.  We introduce a notion of equilibrium that we refer to as trembling-hand-perfect mean field equilibrium (T-MFE), and show existence of T-MFE and a globally convergent computational method in games with strategic complementarities.\vspace{0.05in}\\
(iii) We propose two learning algorithms---model-free and model-based---for learning T-MFE under an unknown model.  The algorithms follow a structure of identifying TQ-values to find a candidate strategy, and then taking only one-step Mckean-Vlasov dynamics to update the mean field distribution.  We show convergence and determine their sample complexity bounds.\vspace{0.05in}\\
(iv) Finally, to the best of our knowledge, ours is the first work that develops a fully online learning scheme that utilizes the large population of agents to concurrently sample the real world.  Our algorithm only needs a one-step Mckean-Vlasov mean field update under each strategy, which automatically happens via agents applying the current strategy for one step.  This obviates the need for a multi-step simulator required by typical RL approaches.  


\paragraph{Related Work:}  There has recently been much interest in the intersection of machine learning and collective behavior, often under the large population regime.  Much of this work focuses on specific classes of stochastic games that possess verifiable properties on information, payoffs and preferences that provide structure to the problem.  In line with this approach is learning in structured MFGs, such as in linear-quadratic, oscillator or potential game settings~\citep{kizilkale2012mean,yin2013learning,cardaliaguet2017learning,carmona2019linear}.  Other approaches provide structure to the problem by considering localized effects, such as local interactions \citep{yang2018mean}, local convergence \citep{mguni2018decentralised}, or a local version of Nash equilibrium \citep{subramanian2019reinforcement}.    While we are not aware of any work that considers learning in games with strategic complementarities, a survey on the algorithmic aspects of multi-agent RL, including a review of existing work in the mean field domain is available in~\cite{zhang2019multi}.

Many of the issues faced in simultaneous learning and decision making in games with large populations are contained in \cite{yang2018mean,guo2019learning,subramanian2019reinforcement}, which are the closest to our work.  \cite{yang2018mean} considers the scenario wherein interactions are local in that each agent is impacted only by the set of neighbors, and so sampling only among them is sufficient to obtain a mean field estimate.   This, however, requires a specific structure in which the Q-function of agents can be decomposed into such local versions.  \cite{subramanian2019reinforcement} proposed a simulator-based (not fully online) policy gradient algorithm for learning the mean field equilibrium in stationary mean field games. However, the convergence guarantees are limited to only a local version of Nash equilibrium.   Finally, \cite{guo2019learning} presents an existence and simulator-based model-free learning algorithm for MFE without structural assumptions on the game.  Instead, there are contraction assumptions imposed on the trajectories of state and action distributions over time, which, however, may not be verifiable for a given game in a straightforward manner.

Our work is distinguished from existing approaches in several ways.  We introduce the concept of a \emph{strategically consistent} approach to learning via the trembling-hand idea, rather than arbitrarily adding a modicum of exploration to best responses as do most existing works.  We provide a structured application scenario to apply this idea in the form of games with strategic complementarities.  In turn, this allows us to explore both model-free and model-based methods to compute and learn optimal trembling-hand perfect MFE, including showing global convergence and determining their sample complexities.  Perhaps most importantly, we exploit the large population setting to obtain samples without the aid of a simulator, which in turn enables learning directly from the real system.

\section{Mean Field Games and Trembling-Hand Perfection}
\label{sec:formulation}
\paragraph{Mean Field Games:}
An $N$-agent stochastic dynamic game is represented as $(\mathcal{S}, \mathcal{A}, P, (r^{i})^{N}_{i}, \gamma),$ where $\mathcal{S}$ and $\mathcal{A}$ are the state and action spaces, respectively, both assumed to be finite. At time $k,$ agent $i$ has state $s^{i}_{k} \in \mathcal{S},$  takes action $a^{i}_{k} \in \mathcal{A},$ and receives a reward $r^{i}(s_{k}, a_{k}).$ Here, $s_{k} = (s^{i}_{k})^{N}_{i=1}$ is the system state and $a_{k} = (a^{i}_{k})^{N}_{i=1}$ is the joint action. The system state evolves according to transition kernel $s_{k+1} \sim P(\cdot | s_{k}, a_{k})$.  Each agent aims at maximizing the infinite horizon cumulative discounted reward $\mathbb{E}[\sum^{\infty}_{k=0} \gamma^{t} r^{i}(s_{k}, a_{k})],$ with discount factor $\gamma \in (0, 1).$ 

Identification of a best response is computationally hard under a Bayesian framework, and a more realistic approach, aligned with a typical agent's computational capabilities is to reduce the information state of each agent to the so-called \emph{mean field} state distribution~\citep{lasry2007mean,tembine2009mean,adlakha2013mean,iyer2014mean,li2016incentivizing,LiXia18}.  The mean field $z_k$ at time $k$ is defined as $z_{k} (s) = \frac{1}{N} \sum^{N}_{i=1} \mathbf{1}\{s^{i}_{k} = s\},$ and is the empirical distribution of the states of all agents.  It represents the agent's belief that the states of all others will be drawn in an i.i.d. manner from $z_k$.  Agent $i$ at time $k$ receives a reward $r^{i}(s^{i}_{k}, a^{i}_{k}, z_{k}),$ and its state evolves according to $s^{i}_{k+1} \sim P(\cdot | s^{i}_{k}, a^{i}_{k}, z_{k})$.  The mean field approximation is accurate under structural assumptions on correlation decay across agent states as the number of agents becomes asymptotically large~\citep{graham1994chaos,iyer2014mean}.

We represent a MFG as $\Gamma = (\mathcal{S}, \mathcal{A}, P, r, \gamma),$ and restrict our attention to stationary MFGs with a homogeneous reward function.  Here, all agents follow the same stationary strategy $\mu :  \mathcal{S} \rightarrow  \mathcal{P}(\mathcal{A}),$ where $\mathcal{P}(\mathcal{A})$ is the probability distribution over the action space, and reward function $r^{i} = r$ for all $i.$ We assume that $|r(s,a,z)| \leq 1.$  The mean field $z_{k}$ evolves following the discrete time McKean-Vlasov equation
\begin{align}
\label{eq:Mckean}
& z_{k+1} = \Phi(z_{k}, \mu),~\text{where,}~ z_{k+1}(s') = \Phi(z_{k}, \mu)(s') \\ & = \sum_{s \in \mathcal{S}} \sum_{a \in \mathcal{A}} z_{k}(s) \mu(s, a) P(s' |s, a, z_{k}).
\end{align}
The value function $V_{\mu, z}$ corresponding to the strategy $\mu$ and the mean field $z$ is defined as $V_{\mu, z}(s) = \mathbb{E}[\sum^{\infty}_{k=0} \gamma^{t} r(s_{k}, a_{k}, z) | s_{0} =s],~\text{with,}~ a_{k} \sim \mu(s_{k}, \cdot), s_{k+1} \sim P(\cdot | s_{k}, a_{k}, z). $
\paragraph{Mean Field Games with Strategic Complementarities (MFG-SC):}
Structural assumptions on the nature of the MFG are needed in order to show existence of an equilibrium, and to identify convergent dynamics.  
Our focus is on a particular structure called \emph{strategic complementarities} that aligns the increase of an agent's strategy with increases the competitors' strategies~\citep{nowak2007stochastic,vives2009strategic,adlakha2013mean}.

We introduce some  concepts before defining MFG-SC.  The partially ordered set $(X,\succeq)$ is called a $\textit{lattice}$ if for all $x,y \in X$, the elements $\sup \{x,y \}$ and $\inf \{x,y \}$ are in $X$.  $X$ is a $\textit{complete lattice}$ if for any  $S \subset X$, both $\sup S$ and $\inf S$ are in X. A function $f : X \rightarrow \mathbb{R}$ is said to be $\textit{supermodular}$ if $f(\sup \{x,x' \}) + f(\inf \{x,x' \}) \ge f(x) + f(x')$ for any  $x,x' \in X$. Given  lattices $X, Y$, a function $f : X \times Y \rightarrow \mathbb{R}$ is said to have $\textit{increasing differences}$ in $x$ and $y$ if for all $x' \succeq x, y' \succeq y$,  $f(x',y') - f(x',y) \ge f(x,y')-f(x,y)$. A correspondence $T: X \rightarrow Y$ is $\textit{non-decreasing}$ if $x' \succeq x, y \in T(x),$ and $y' \in T(x')$ implies that $\sup \{y,y' \} \in T(x')$ and $\inf \{y,y' \} \in T(x)$.

For probability distributions $p, p' \in \mathcal{P}(X)$, we say $p$ \textit{stochastically dominates} $p'$, denoted as $p \succeq_{SD} p'$,   if $\sum_{x \in X} f(x) p(x)  \geq  \sum_{x \in X} f(x) p'(x)$ for any bounded non-decreasing function $f$.  The conditional distribution $p(\cdot | y)$ is  stochastically non-decreasing in $y$ if for all $y' \succeq y$, we have $p(\cdot | y') \succeq_{SD} p(\cdot | y)$. Finally, the conditional distribution $p(\cdot | y, z)$ has \textit{stochastically increasing differences} in $y$ and $z$ if $\sum_{x} f(x) p(x|y,z)$ has increasing differences in $y$ and $z$ for any bounded non-decreasing function $f$.

We now give the formal definition of mean field games with strategic complementarities \citep{adlakha2013mean}. 
\begin{definition}[MFG-SC]\label{def:mfg-sc}
Let $\Gamma$ be a stationary mean field game. We say that $\Gamma$ is a mean field game with strategic complementaries if it has the  following properties:\\
(i)  Reward function:  $r(s, a, z)$ is  non-decreasing in $s$,  supermodular in $(s,a)$, and has increasing differences in $(s,a)$ and $z$. Also, $\max_{a } r(s,a,z)$ is non-decreasing in $s$ for all fixed $z$. \\
(ii) Transition probability:  $P(\cdot | s,a,z)$ is stochastically supermodular in $(s,a)$, has stochastically increasing differences in $(s,a)$ and $z$, and  is stochastically non-decreasing in each of $s,a,$ and $z$.
\end{definition}

\paragraph{MFG-SC Example - Infection Spread:}  We assume a large but fixed population of agents (computers or humans). At any time step, an agent may leave the system (network or town) with probability $\zeta,$ and is immediately replaced with a new agent.   The state of an agent $s\in \mathbb{Z}^+$ is its health level, and the agent can take action $a\in \{1, 2, 3,...|\mathcal{A}|\}$ (installing security patches, wearing a mask etc.) to stay healthy. The susceptibility of each agent, $p(s)$, is a decreasing function in state (higher health implies lower susceptibility).  At each time step, the agent interacts with the ensemble of agents who have a mean field state distribution $z.$  We define infection intensity $i_{z}=c_f p(\sum_{s \in \mathcal{S}}sz(s)),$ which can be interpreted as the probability of getting infected via interaction with the population and $c_f$ is the infection intensity constant.

Given the current state-action pair $(s,a),$ the next state $s'$ is given by \[s' = (s+a-w_1)_{+}\mathbbm{1}\{E_1\} + (s+a)\mathbbm{1}\{E_2\} + w_2 \mathbbm{1}\{E_3\},\] where $(x)_+ = \max\{0, x\}$ (state is non-negative),  $w_1,w_2 \in \mathbb{Z}^+$ are realizations of non-negative random variables, and $E_{i}$s are mutually exclusive events with probabilities $i_{z} (1-\zeta), (1-i_{z})(1-\zeta),$ and $\zeta,$ respectively. {Events $E_1$ and $E_2$ correspond to the agent remaining in the system, and being infected (health may deteriorate) or not infected, respectively.}  $E_3$ is the event that the agent leaves and is replaced with an agent with random state (regeneration).  An agent receives a reward that depends on its own immunity $1-p(s)$ as well as that of the population (i.e., system value increases with immunity), but pays for its action. Hence,  \[r(s,a,z) = \delta_1 (1 -p(s)) + \delta_2 \sum_{s \in \mathcal{S}}z(s) (1 -p(s)) - \delta_3 a,\] where $\delta_{i}$s are positive constants. 

It is easy to verify that the model has non-decreasing differences in the transition matrix.  Also, it encourages crowd seeking behavior, where the mean field positively affects rewards obtained.  Showing that it has strategic complementarities is straightforward. We validate our algorithms via simulations on this model in Section \ref{sec:simulations}.


\paragraph{MFG-SC Example - Amazon Mechanical Turk (MTurk) and other Gig Economy Marketplaces:}  
We consider Amazon Mechanical Turk (MTurk), a crowd sourcing market wherein human workers (called Turkers) are recruited to perform so-called Human Intelligence Tasks (HITs).  There is a natural alignment of effort employed by Turkers in MTurk, since higher efforts translate into more HITs done right, which results in a higher quality of work distribution, which results in firms willing to spend more on HITs.   
Note that free riding is difficult, since poor quality work results in payments being withheld and reputation loss.  This notion of incentive alignment applies to essentially all Gig economy marketplaces such as Uber and Airbnb---the reputation of the agent directly enhances its reward, while the reputation of the marketplace as a whole (i.e., its mean field) draws customers willing to pay into the system, and so enhances the reward of the agent. A more detailed description and numerical simulations  are provided in the supplementary material. 

\paragraph{Trembling-Hand-Perfect Mean Field Equilibrium:}
The notion of trembling-hand-perfection is a means of refining the Nash equilibrium concept to account for the fact that equilibria that naturally occur are often those that are optimal when a known amount of randomness is introduced into the strategies employed to ensure that they are totally mixed, i.e., all actions will be played with some (however small) probability~\citep{bielefeld1988reexamination}.  Thus, the agent is restricted to only playing such  mixed (randomized) strategies, but  maintains strategic consistency in that it accounts for the probability of playing each action while calculating the expected payoff of such a totally mixed strategy. 

Formalizing the above thoughts, we denote the set of trembling-hand strategies as $\Pi^{\epsilon}$. A trembling-hand strategy $\mu \in \Pi^{\epsilon}$ is a mapping $\mu : \mathcal{S} \rightarrow \mathcal{P}^{\epsilon}(\mathcal{A}),$ where $\mathcal{P}^{\epsilon}(\mathcal{A})$ is the set of $\epsilon$-randomized probability vectors over $\mathcal{A}$.  Any probability vector in $\mathcal{P}^{\epsilon}(\mathcal{A})$ has the value $(1-\epsilon)$ for one element and the value $\epsilon/(|\mathcal{A}|-1)$ for all the other  elements. Thus, any $\mu \in \Pi^{\epsilon}$ has the following form: $\mu(s, a) = (1 - \epsilon)$ for $a = a_{s}$ for some $a_{s}$ and $\mu(s, a) =  \epsilon/(|\mathcal{A}|-1)$ for all other $a \in \mathcal{A}$. In the standard reinforcement learning parlance, $\Pi^{\epsilon}$ is essentially the set of all $\epsilon$-greedy policies.

The main difference between a trembling-hand strategy and an $\epsilon$-greedy policy lies in the value function.  Recall that under the $\epsilon$-greedy idea, the agent computes the pure (deterministic) best response policy, and then arbitrarily adds randomization.  However, under the strategic game setting, choosing a strategy that could result in an arbitrary loss of value is impermissible.   Rather, the agent must compute the best trembling-hand strategy, i.e., it must account for the impact on value of the $\epsilon$ randomness.  Formally, we first define the optimal trembling-hand value function $V^{*}_{z}$ and the optimal trembling hand strategy $\mu^{*}_{z}$ corresponding to the mean field $z$ as
\vspace{0.0cm}
\begin{align}
    V^{*}_{z} &= \max_{\mu \in \Pi^{\epsilon}} V_{\mu, z},\\
    \mu^{*}_{z} &\in \Psi(z),~~\text{where}~~\Psi(z) = \arg \max_{\mu \in \Pi^{\epsilon}} V_{\mu, z}.
\end{align}
We define trembling-hand-prefect mean field equilibrium (T-MFE) in terms of a trembling-hand strategy $\mu^*$ and a mean field distribution $z^*$ that must jointly satisfy, (i) \emph{optimality}---the strategy $\mu^*$ must be superior to all other strategies, given the belief $z^*,$ and (ii)  \emph{consistency}---given a candidate mean field distribution $z^*,$ the strategy $\mu^*$ must regenerate $z^*$ under the McKean-Vlasov dynamics~(\ref{eq:Mckean}).
\begin{definition}[T-MFE]
Let $\Gamma$ be a stationary mean field game. A trembling-hand-perfect  mean field equilibrium $(\mu^{*}, z^{*})$ of ~$\Gamma$ is a strategy $\mu^{*} \in \Pi^{\epsilon}$ and a mean field $z^{*}$ such that,    $\text{(optimality condition)}~~ \mu^{*}  \in \Psi(z^{*}) ,~\text{and}~ \text{(consistency condition)}~~z^{*} = \Phi(z^{*}, \mu^{*})$
\end{definition}
\section{Existence and Computation of T-MFE}\label{sec:existence-comp}


\paragraph{Existence of T-MFE:}
We first introduce a method to compute the optimal trembling-hand value function $V^{*}_{z}$ for a given mean field $z$.  Note that classical value iteration for any given finite MDP will converge under deterministic policies (pure strategies)---something that is not possible under our restriction to trembling hand strategies $\Pi^{\epsilon},$ which only allows totally mixed strategies.   We overcome this issue by using a generalized value iteration approach~\citep{szepesvari1996generalized}.  
%
%
Here, rather than the value function, we compute the Q-value function, which for a strategy $\mu$ and a given mean field $z$ is defined as $Q_{\mu, z}(s, a) = \mathbb{E}[\sum^{\infty}_{t=0} \gamma^{t} r(s_{t}, a_{t}, z) | s_{0} =s, a_{0} = a],~\text{with,}~ a_{t} \sim \mu(s_{t}, \cdot), s_{t+1} \sim P(\cdot | s_{t}, a_{t}, z).$
The optimal trembling-hand Q-value function (TQ-value function) for a given mean field $z$ is then defined as $Q^{*}_{z} = \max_{\mu \in \Pi^{\epsilon}}  Q_{\mu, z}$. The optimal trembling-hand strategy $\mu^{*}_{z}$ for a given mean field $z$ can then be computed as $\mu^{*}_{z} = \pi^{\epsilon}_{Q^{*}_{z}}.$

For any given $Q$-value function, define the trembling-hand strategy  $\pi^{\epsilon}_{Q}$ and the  function  $G(Q)$ as
\begin{align*}
\pi^{\epsilon}_{Q}(s, a) & = \left \{ \begin{array}{ll}
         (1 - \epsilon)& \text{for}~a =   \arg \max_{b} Q(s, b)  \\
          {\epsilon}/({|\mathcal{A}| - 1})& \text{for}~a \neq \arg \max_{b} Q(s, b)
    \end{array} \right., \\
    G(Q)(s) & = \sum_{a \in \mathcal{A}} \pi^{\epsilon}_{Q}(s, a) ~ Q(s, a).
\end{align*}
Note that $\pi^{\epsilon}_{Q}$ is the usual $\epsilon$-greedy policy with respect to $Q$. Using the above notation, we  define the TQ-value operator  $F_{z}$ for a given mean field $z$ as
\begin{align}
\label{eq:F-operator}
    F_{z}(Q)(s, a) = r(s, a, z) + \gamma \sum_{s'} P(s' | s, a, z) ~G(Q)(s'). 
\end{align}
The TQ-value operator $F_{z}$ has properties similar to the standard Bellman operator. In particular, we show below that $ F_{z}$ is a contraction with $Q^{*}_{z}$ as its unique fixed point. 
\begin{proposition}
\label{prop:F-mapping}
(i) $F_{z}$ is a contraction mapping in sup norm  for all $z \in \mathcal{P}(\mathcal{S})$. More precisely,   $\|F_{z}(Q_{1}) -  F_{z}(Q_{2})\|_{\infty} \leq \gamma  \|Q_{1} - Q_{2} \|_{\infty}$ for any $Q_{1}, Q_{2},$ and for all $z \in \mathcal{P}(\mathcal{S})$. \\
(ii) The optimal trembling-hand Q-value function $Q^{*}_{z}$~ for a given mean field $z$ is the unique fixed point of $F_{z}$, i.e., $F_{z}(Q^{*}_{z}) = Q^{*}_{z}$.
\end{proposition}
The contraction property of $F_{z}$ implies that the iteration  $Q_{m+1, z} = F_{z}(Q_{m,z})$  will converge to the unique fixed point of $F_{z}$, i.e., $Q_{m,z} \rightarrow Q^{*}_{z}$.  We call this procedure as TQ-value iteration. 



From the above result, we can compute the optimal trembling-hand strategy for a given mean field $z.$   However, it is not clear if there exists a T-MFE $(z^{*}, \mu^{*})$ that simultaneously satisfies the optimality condition and consistency condition. We answer this question affirmatively below.
\begin{theorem}
\label{thm:tmfe-existence}
Let $\Gamma$ be a stationary mean field game with strategic complementarities. Then, there exists a trembling-hand-perfect mean field equilibrium for $\Gamma.$
\end{theorem}
The proof  follows  from the monotonicity properties of MFG-SC.  Thus, given this game structure, no additional conditions are needed to show existence.


\paragraph{Computing T-MFE:}
Given that a T-MFE exists, the next goal is to devise an algorithm to compute a T-MFE.  A natural approach is to use a form of best-response dynamics as follows. Given the candidate mean field $z_{k}$, the trembling best-response strategy $\mu_{k}$ can be computed as $\mu_{k} \in \Psi(z_{k})$. While the typical approach would be to then compute the stationary distribution under $\mu_{k},$ we simply update the next mean field $z_{k+1}$ by using just one-step Mckean-Vlasov dynamics as $z_{k+1} = \Phi(z_{k}, \mu_{k}),$ and the cycle continues. 
While the approach is intuitive and reminiscent of the best-response dynamics proposed in \cite{adlakha2013mean}, it is not clear that it will converge to any equilibrium.  
We show that in mean field games with strategic complementarities, such a trembling best-response (T-BR) process converges to a T-MFE.  Our computation algorithm, which we call the T-BR algorithm, is presented in Algorithm~\ref{algo:T-MRD}. 
\begin{algorithm}
\caption{T-BR Algorithm}
\label{algo:T-MRD}
\begin{algorithmic}[1]
\STATE Initialization: Initial mean field $z_{0}$ 
\FOR {$k=0,1,2,..$}
\STATE For the mean field $z_{k}$, compute the optimal TQ-Value function $Q^{*}_{z_{k}}$ using the TQ-value iteration $Q_{m+1,z_{k}} = F(Q_{m,z_{k}})$
\STATE Compute the strategy $\mu_{k} = \Psi(z_{k})$ as the trembling-hand strategy w.r.t $Q^{*}_{z_{k}}$, i.e, 
$\mu_{k} =  \pi^{\epsilon}_{Q^{*}_{z_{k}}}$ 
\STATE Compute the next mean field  $z_{k+1} = \Phi(z_{k}, \mu_{k})$
\ENDFOR 
\end{algorithmic}
\end{algorithm}

\begin{theorem}
\label{thm:T-MRD}
Let $\Gamma$ be a stationary mean field game with strategic complementarities. Let $\{z_{k}\}$ and  $\{\mu_{k}\}$  be the sequences of mean fields and strategies  generated according to Algorithm \ref{algo:T-MRD}. Then  $\mu_{k} \rightarrow \mu^*$   and $z_{k} \rightarrow z^*$  as $k \rightarrow \infty$ where  $(\mu^*,z^*)$ constitutes a T-MFE of ~$\Gamma.$ 
\end{theorem}

\begin{remark}
\label{remark:k0}
Consider the sequence of mean fields $\{{z}_{k}\}$ generated by the T-BR algorithm.  According to Theorem \ref{thm:T-MRD}, there exists a finite $k_{0} = k_{0}(\bar{\epsilon})$ such that  $\|{z}_{k} - z^{*}\| \leq \bar{\epsilon}$~ for all $k \geq k_{0}$, where $z^{*}$ is the T-MFE. A precise characterization of $k_{0}$ is difficult because the convergence of the T-BR algorithm is based on monotonicity properties of MFG-SC, rather than on contraction arguments. 
\end{remark}




\section{Learning T-MFE}
\label{sec:learning-TMFE-offline}
We address the problem of learning T-MFE when the model is unknown. In this section, we assume the availability of a simulator, which, given the current state $s$, current action $a$ and current mean field $z$, can generate the next state $s' \sim P(\cdot | s, a, z)$. We discuss how to learn directly from real-world samples without a simulator in the next section.  We also assume that the reward function is known, as is common in the literature.  We now introduce two reinforcement learning algorithms---a model-free algorithm and a model-based algorithm---for learning  T-MFE.

\subsection{Model-Free TMFQ-Learning Algorithm for learning T-MFE}\label{sec:QL-conv-TMFE}
We first describe the TMFQ-learning algorithm, which builds on the T-BR algorithm. Recall that the T-BR algorithm uses knowledge of the model in two locations.  The first is that at each step $k,$ for a given mean field $z_{k}$, the  optimal  TQ-value function $Q^{*}_{z_{k}}$ is computed using TQ-value iteration. 
In the learning approach, we use the generalized Q-learning framework \citep{szepesvari1996generalized} as the basis for the model-free TQ-learning algorithm as follows:
\begin{align}
\label{eq:TQ-update}
\begin{split}
    & Q_{t+1, z_k}(s_{t}, a_{t}) = (1 - \alpha_{t})  Q_{t, z_k}(s_{t}, a_{t}) \\
    & \quad \quad \quad + \alpha_{t} (r(s_{t}, a_{t}, z_{k}) + \gamma G(Q_{t, z_k})(s_{t+1}))
\end{split}
\end{align}
where $\alpha_{t}$ is the appropriate learning rate. Here, the state sequence $\{s_{t}\}$ is generated using the simulator by fixing the mean field $z_{k},$ i.e., 
$s_{t+1} \sim P(\cdot | s_{t}, a_{t}, z_{k}), \forall t$. Using the properties of the generalized Q-learning formulation \citep{szepesvari1996generalized}, it can be shown that $Q_{t, z_k} \rightarrow Q^{*}_{z_{k}}$ as $t \rightarrow \infty$. 

The second location where the model is needed in T-BR is for the one-step McKean-Vlasov update.  In the learning approach, given the  mean field $z_{k}$ and the strategy $\mu_{k}$, the next mean field $z_{k+1}$ can be estimated to a desired accuracy using the simulator.  The precise numerical approach is presented as the {Next-MF} scheme described in  Algorithm \ref{algo:next-MF} {(supplementary material).}  We can now combine these steps to obtain the TMFQ-learning algorithm presented in Algorithm~\ref{algo:TQL}. 
\begin{algorithm} 
\label{algo:TQL}
	\caption{TMFQ-Learning Algorithm for T-MFE}
	\begin{algorithmic}[1]
		\STATE Initialization:  Initial mean field ${z}_0$
		\FOR {$k=0, 1, 2, \ldots $}
		    \STATE Initialize time step $t \leftarrow 0$. Initialize $s_{0}, Q_{0,z_k}$
		    \REPEAT
		        \STATE Take action $a_{t} \sim \pi^{\epsilon}_{Q_{t,z_k}}(s_{t})$, observe  reward and  the next state $s_{t+1} \sim P(\cdot | s_{t}, a_{t}, z_{k})$
		        \STATE Update $Q_{t,z_k}$ according to TQ-learning \eqref{eq:TQ-update} 
		        \STATE $t \leftarrow t+1$
		    \UNTIL{$\norm{Q_{t,z_k}-Q_{t-1,z_k}}_{\infty} < \epsilon_{1}$}
		    \STATE Let $Q_{z_{k}} = Q_{t,z_k}$ and the strategy $\mu_{k} = \pi^{\epsilon}_{Q_{z_k}}$
		    \STATE $z_{k+1} =${Next-MF}($z_k,\mu_{k}$)
		\ENDFOR
	\end{algorithmic}
\end{algorithm}

We first present an asymptotic convergence result of TMFQ-learning based on a perfect accuracy assumption on the TQ-learning and Next-MF steps, i.e., they are run to convergence. A complex, but more accurate analysis using two-timescale stochastic approximation is also possible.  Instead, we will remove this assumption, and provide a PAC-type result further below.
\begin{theorem}
\label{thm:T-MFE-convergence}
Let $\Gamma$ be a stationary mean field game with strategic complementarities. Let $\{z_{k}\}$ and $\{\mu_{k}\}$ be the sequences of the mean fields and policies generated by Algorithm \ref{algo:TQL}. Then ${\mu}_{k} \rightarrow {\mu}^*$  and $z_{k} \rightarrow z^{*}$  as $k \rightarrow \infty$ where $(\mu^*, z^{*})$ is a T-MFE of ~$\Gamma$. 
\end{theorem}

In a practical implementation, we may only run the TQ-learning step and the Next-MF step for a finite number of iterations. Hence, we develop a sample complexity bound under which TQ-learning and Next-MF provide an appropriate accuracy. We desire to compare TMFQ-learning with T-BR after $k_{0}$ iterations, where $k_0$ is the number of iterations of T-BR that yields a mean field that is $\bar{\epsilon}$-close to the T-MFE $z^{*}$. We make some necessary assumptions that are required for such a characterization.




\begin{assumption}
\label{assump:Lips-R-P}
(i) There exists  $C_{1} > 0$ such that  $\norm{ r(\cdot, \cdot, z) - r(\cdot, \cdot, \tilde{z}) }_{1} \le C_{1}\norm{z - \tilde{z}}_{1}$, for all $z, \tilde{z} \in \mathcal{P}(\mathcal{S})$.\\
(ii) There exists a $C_{2} > 0$ such that  $\norm{ P(\cdot | \cdot, \cdot, z) - P(\cdot | \cdot, \cdot, \tilde{z}) }_{1} \le C_{2}\norm{z - \tilde{z}}_{1}$, for all $z, \tilde{z} \in \mathcal{P}(\mathcal{S})$. \\
(iii) Let $P_{\mu,z}(s'|s) = \sum_{a} \mu(s, a) P(s'|s,a,z)$. Let $\mu_{1}, \mu_{2}$ be the trembling-hand policies corresponding to $Q_{1}, Q_{2}$, i.e., $\mu_{1} =\pi^{\epsilon}_{Q_{1}}, \mu_{2} = \pi^{\epsilon}_{Q_{2}}$. Then there exists a $C_{3} > 0$ such that $\norm{P_{z,\mu_1}-P_{z,\mu_2}}_{1} \leq C_3 \norm{Q_1-Q_2}_{\infty}$ for all $z \in \mathcal{P}(\mathcal{S})$ and for any given $Q_{1}, Q_{2}$. 
\end{assumption}
Assumption \ref{assump:Lips-R-P}.(i) and \ref{assump:Lips-R-P}.(ii) indicate that the reward function and transition kernel are Lipschitz with respect to the mean field, while Assumption \ref{assump:Lips-R-P}.(iii) indicates that the distance between the Markov chains induced by two policies on the same transition kernel are upper bounded by a constant times the distance between their respective Q-functions.  We then have the following.

\begin{theorem} 
\label{thm:TQL-SC}
Let Assumption \ref{assump:Lips-R-P} hold. For any $0 \leq \bar{\epsilon}, \bar{\delta} < 1,$ let $k_{0} = k_{0}(\bar{\epsilon})$. In Algorithm \ref{algo:TQL}, for each $k \leq k_{0}$, assume that TQ-learning (according to \eqref{eq:TQ-update}) update is performed $T_{0}$ number of times where $T_{0}$ is given as
\begin{align}
\label{eq:asynch-Q-sample-1}
    T_{0} & = O \Bigg( \Bigg(\frac{B^{2} L^{1+3w} V^{2}_{\max}}{\beta^{2} \bar{\epsilon}^{2}} \ln \big( \frac{2 B k_{0}|\mathcal{S}||\mathcal{A}|V_{\max}}{\bar{\delta} \beta \bar{\epsilon}} \big)   \Bigg)^{\frac{1}{w}} \nonumber \\ & + \Bigg( \frac{L}{\beta} \ln(\frac{B V_{\max}}{\bar{\epsilon}}) \Bigg)^{\frac{1}{1-w}} \Bigg),
\end{align}
where $B = (1 + C_{2} + C_{3} D)^{k_{0}+1} (C_{3}+1)$,  $D = (C_{1} + \gamma C_{2})/(1-\gamma)$, $V_{\max} = 1/(1-\gamma), \beta = (1-\gamma)/2$,  $L$ is an upper bound on the covering time \footnote{Covering time of a state-action pair sequence is the number of steps needed to visit all state-action pairs starting from any arbitrary state-action pair \cite{even2003learning}.}, and $w \in (1/2, 1)$. Then,
\begin{align*}
    \mathbb{P}(\|z_{k_{0}} - z^{*} \| \leq 2 \bar{\epsilon}) \geq (1 - \bar{\delta}). 
\end{align*}
\end{theorem}

We may also eliminate the dependence of the constant term $B$ on $k_0$ under a contraction assumption on the McKean-Vlasov dynamics $\Phi$ (for instance, following conditions similar to \cite{borkar2013asymptotics}).  
\begin{assumption}
\label{assump:MV-strong}
Let $Q_{1}, Q_{2}$ be two arbitrary Q-value functions and let $\mu_{1} = \pi^{\epsilon}_{Q_{1}}, \mu_{2} = \pi^{\epsilon}_{Q_{2}}$. Let $z_{1}, z_{2}$ be two arbitrary mean fields. Then there exists positive constants $C_{4}$ and $C_{5}$ such that  $ \| \Phi(z_{1}, \mu_{1}) -  \Phi(z_{2}, \mu_{2}) \|_{1} \leq C_{4} \|z_{1} - z_{2} \|_{1} + C_{5} \|Q_{1} - Q_{2} \|_{\infty}.$ Also assume that $(C_{4} + C_{5} D) < 1,$ where $D = (C_{1} + \gamma C_{2})/(1-\gamma)$. 
\end{assumption}
\begin{corollary}
\label{cor:cor-1}
Let Assumption \ref{assump:Lips-R-P} and Assumption \ref{assump:MV-strong} hold. Then, the we obtain the bound on $T_{0}$ as in \eqref{eq:asynch-Q-sample-1} with $B = (C_{5}+1)/(1 - (C_{4} + C_{5} D)),$ which does not depend on  $k_{0}$. 
\end{corollary}
\subsection{Generative Model-Based Reinforcement Learning for T-MFE}
A model-based variant of the T-BR algorithm is also straightforward to construct.  We note that non-stationarity in our system (the model changes at each step) implies that there is no single model.  Hence, our approach follows the generation of a new model each time that the mean field evolves under a T-BR-like procedure.  Full details are presented in Appendix~\ref{sec:MBRL}.

\section{Online Learning of T-MFE}
\label{sec:learning-TMFE-online}
\begin{figure*}[hbt!]
\begin{minipage}{.32\linewidth}
\centering
\includegraphics[width=1\linewidth]{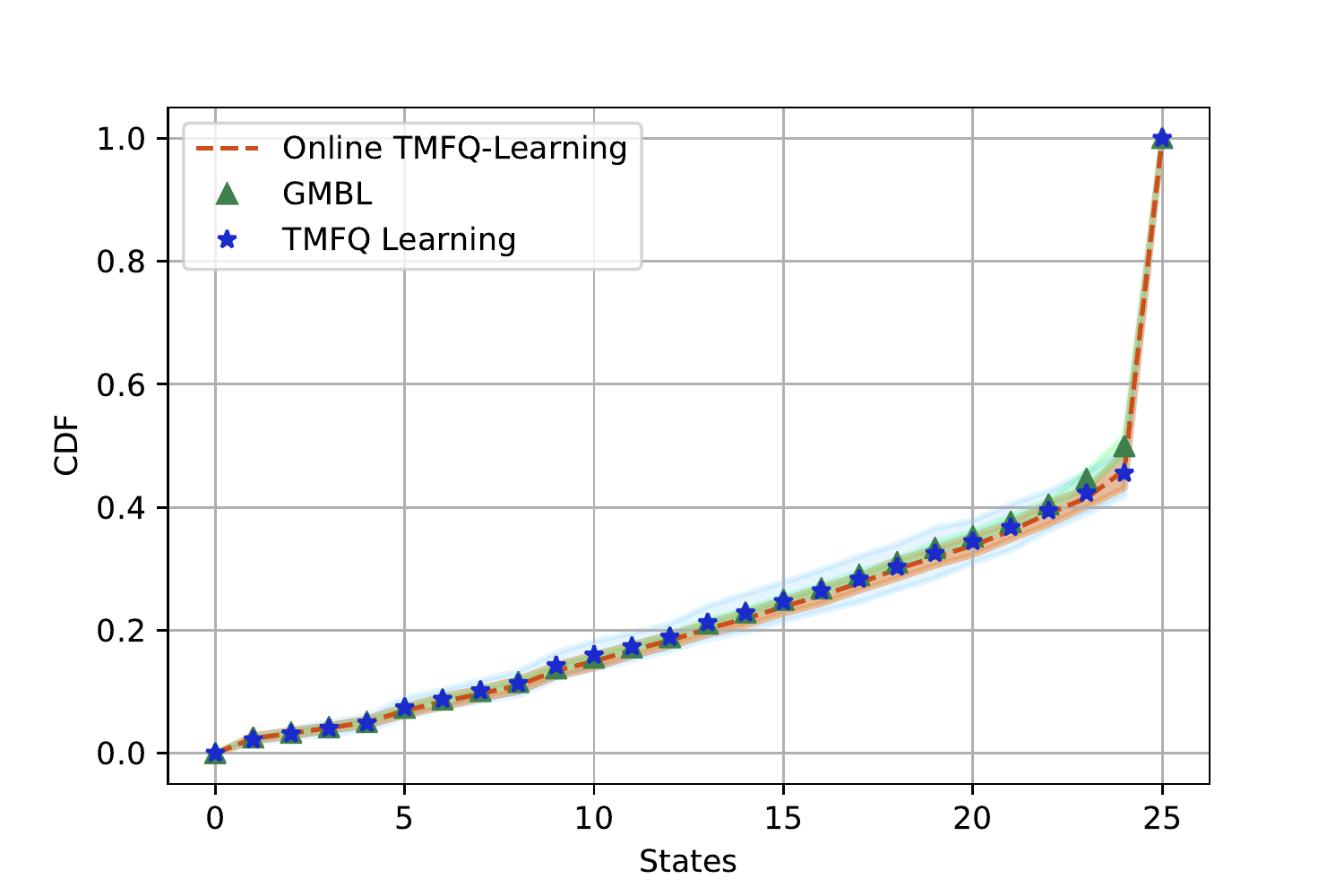}
\caption{TMFE with $c_f=0.1$ }
\label{fig:CDF_Algos}
\end{minipage}\hfill
\begin{minipage}{.32\linewidth}
\centering
\includegraphics[width=1\linewidth]{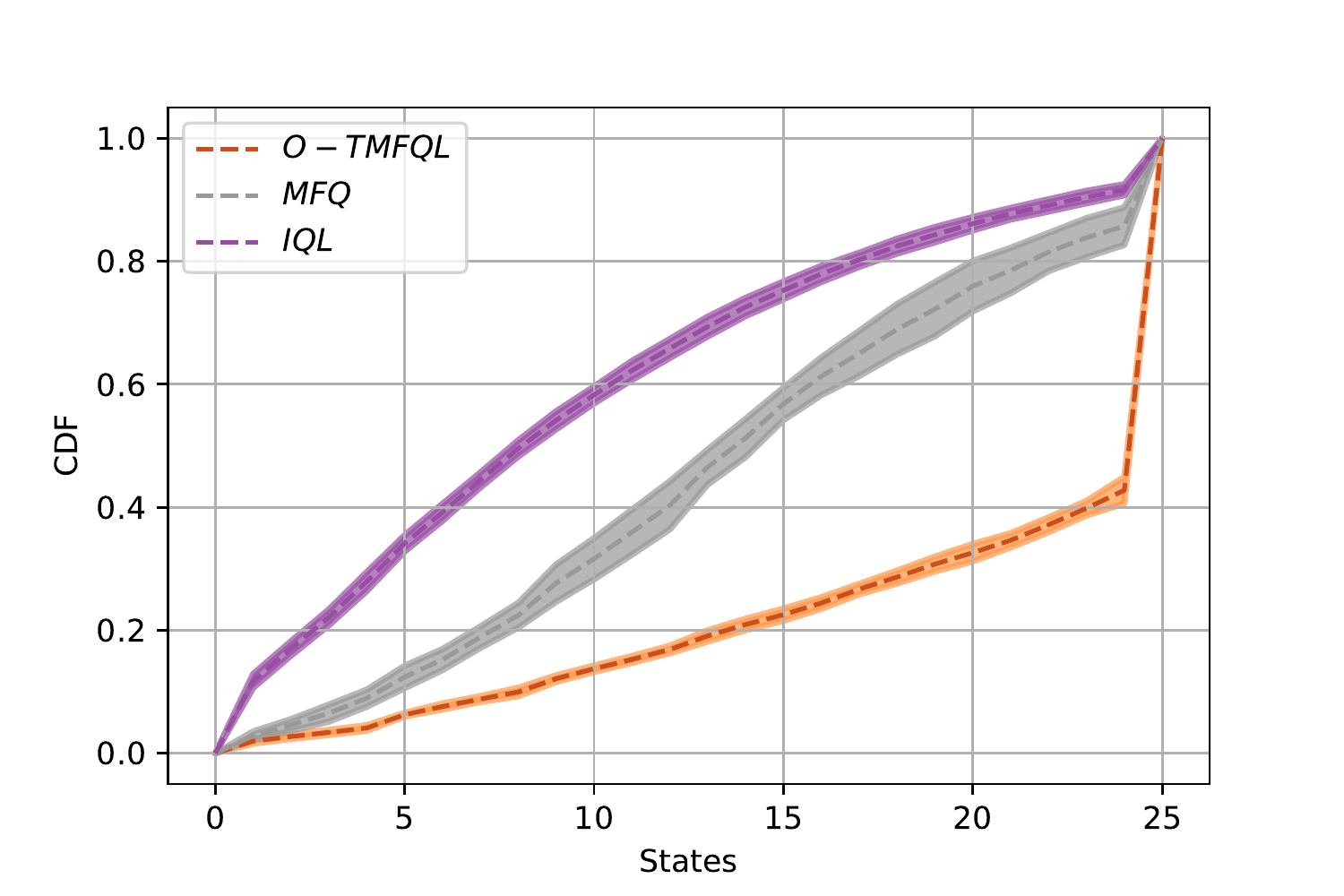}
\caption{CDF of states for different algorithms: $c_f=0.05$ }
\label{fig:Algos_CDF}
\end{minipage}\hfill
\begin{minipage}{.32\linewidth}
\centering
\includegraphics[width=1\linewidth]{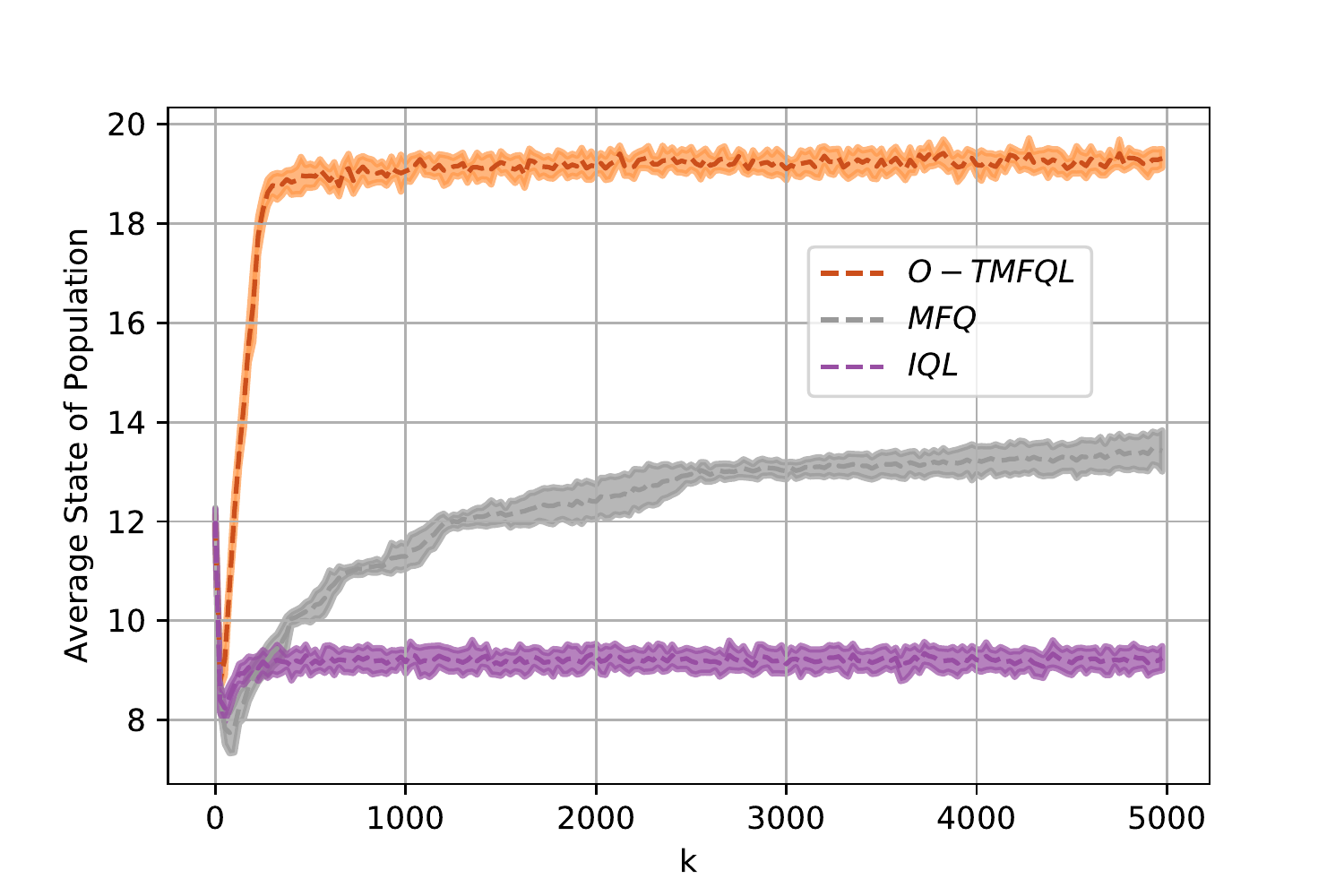}
\caption{Mean states for different algorithms: $c_f=0.05$}
\label{fig:Algos_mean_state}
\end{minipage}\hfill
\begin{minipage}{.32\linewidth}
\centering
\includegraphics[width=1\linewidth]{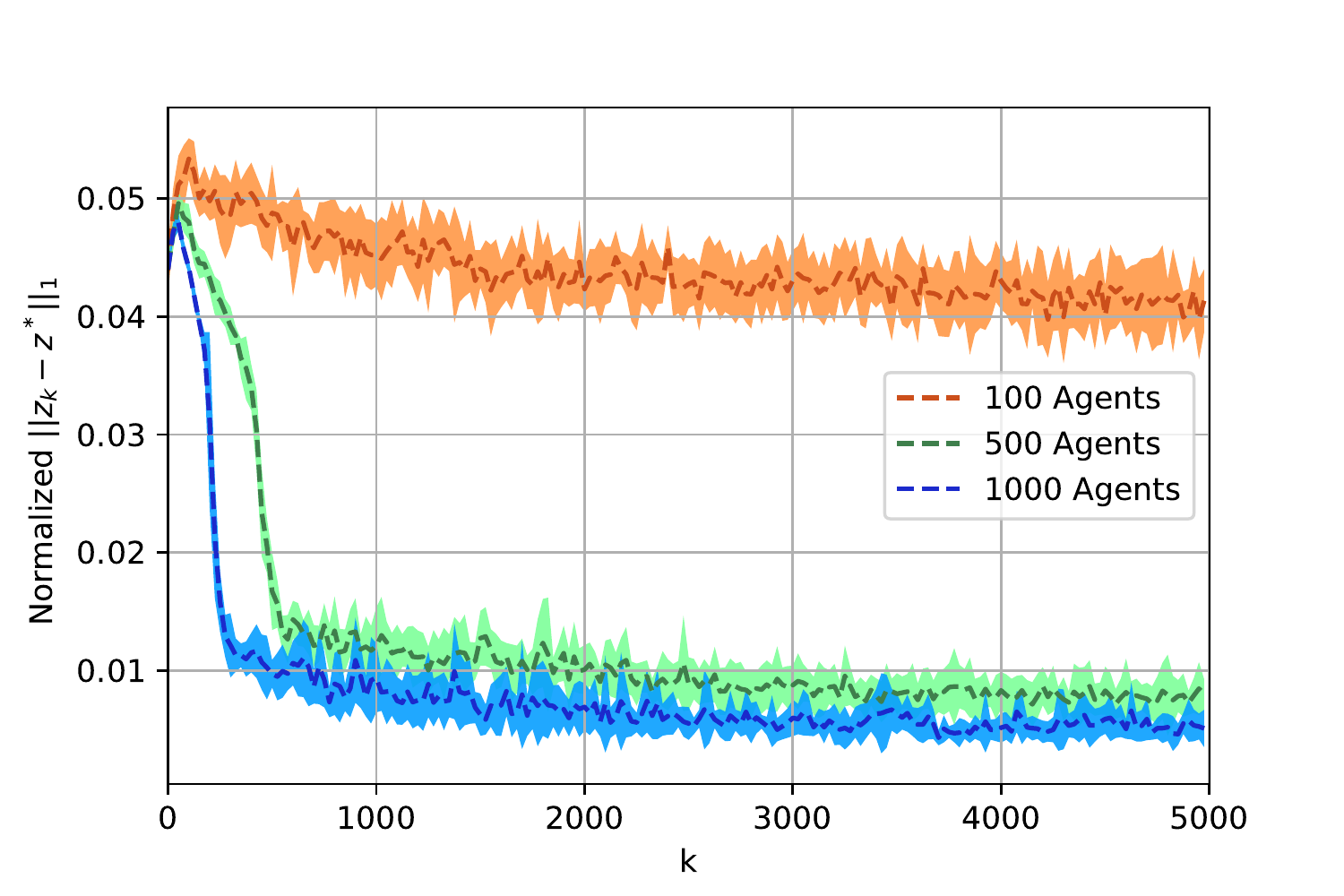}
\caption{Convergence of O-TMFQ-Learning: $c_f=0.05.$ }
\label{fig:z_converge}
\end{minipage}\hfill
\begin{minipage}{.32\textwidth}
\centering
\includegraphics[width=1\linewidth]{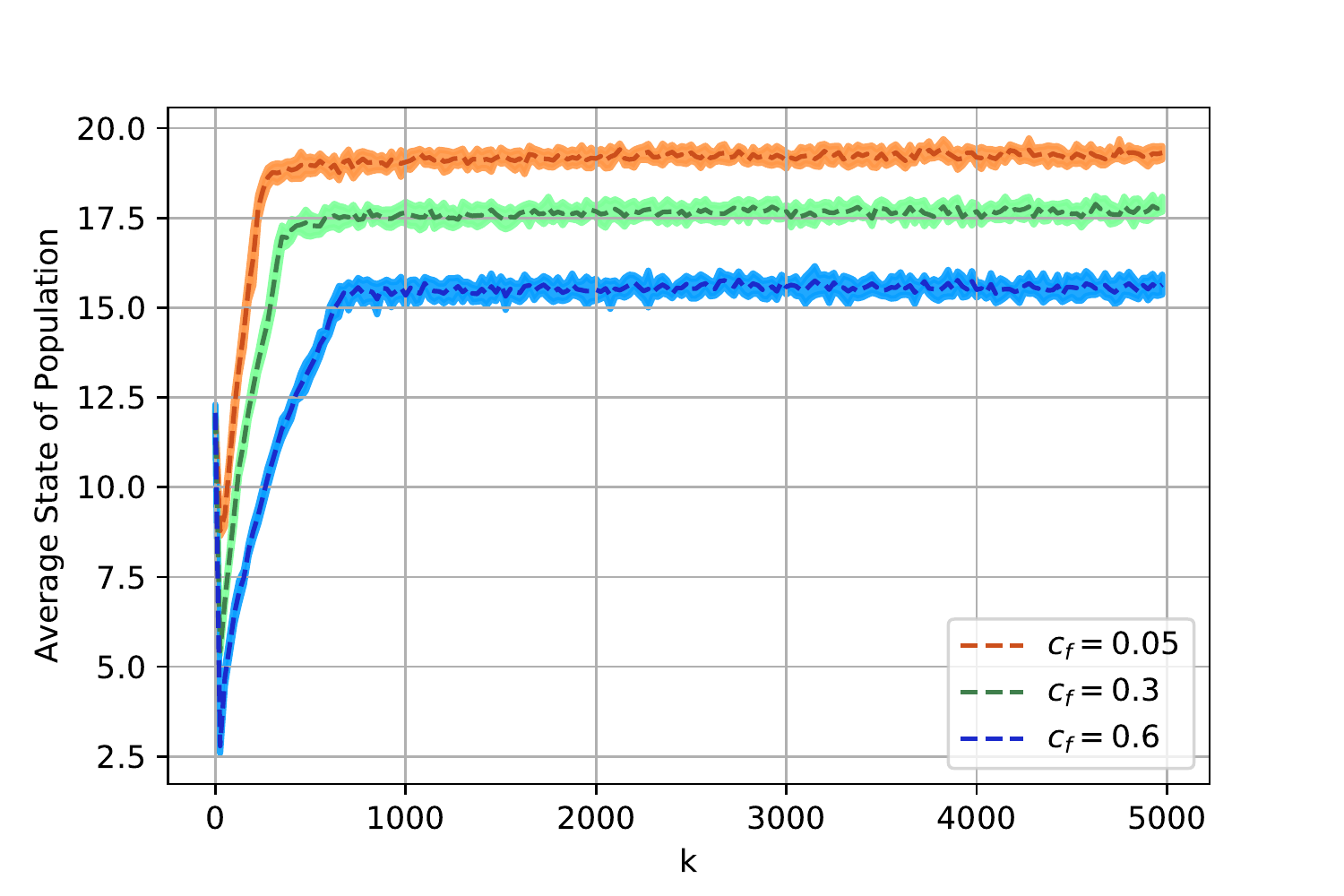}
\caption{Mean state of O-TMFQ-Learning over time}
\label{fig:mean_state}
\end{minipage}\hfill
\begin{minipage}{.32\linewidth}
\centering
\includegraphics[width=1\linewidth]{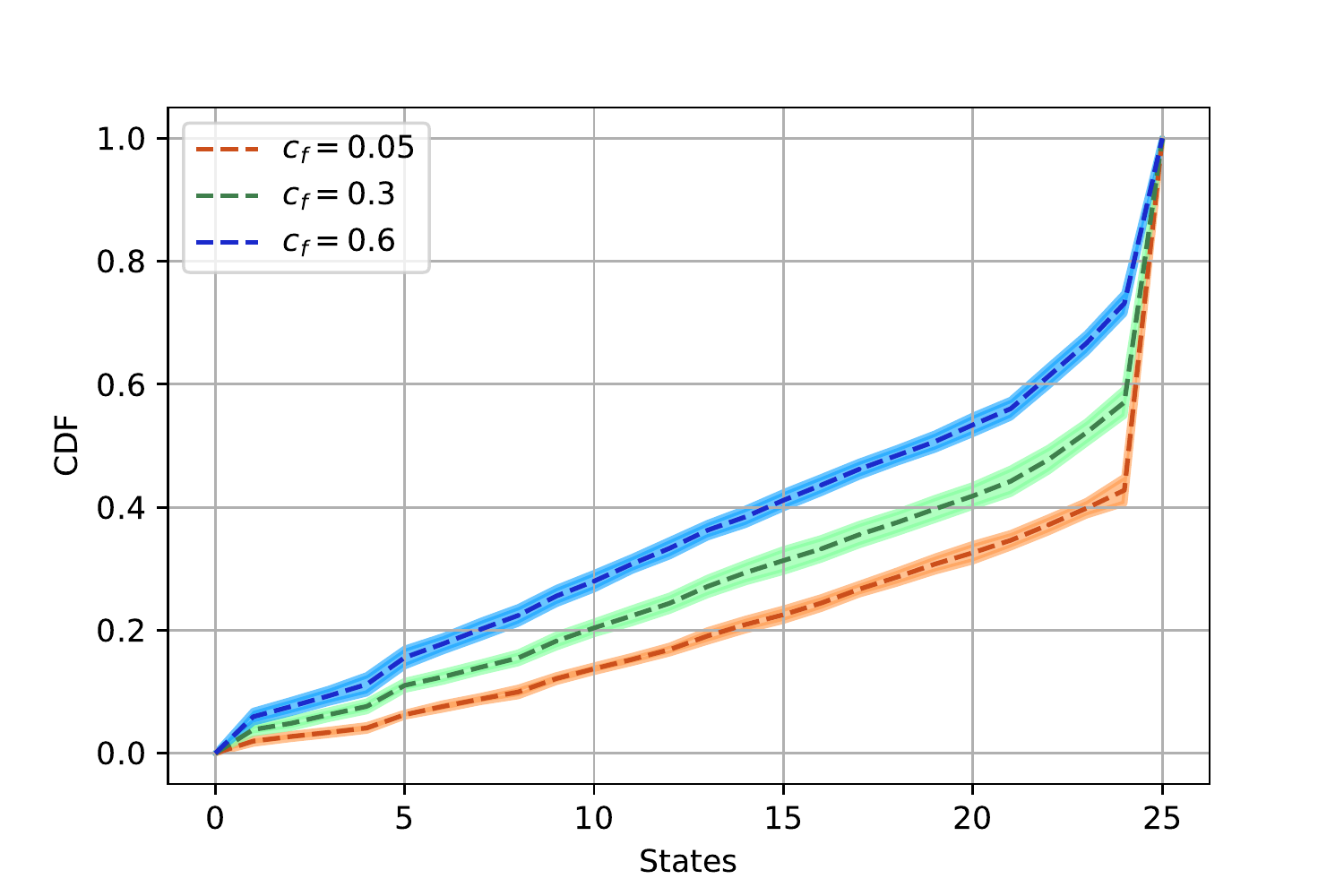}
\caption{CDF of O-TMFQ-Learning}
\label{fig:CDF_OL}
\end{minipage}\hfill
\end{figure*}

The availability of a large number of agents that explore via trembling hand strategies suggests that we can do away with a system simulator via an online algorithm that simply aggregates these concurrently generated samples and computes a new strategy that is then pushed to all agents. Typically, the assumption in such large population scenarios is that the system size is fixed, but any player may leave the system at time step $t$ with probability $\zeta \in (0,1),$ and be immediately replaced by a new agent with random state referred to as a \emph{regeneration event.}

A generic RL approach would require that given the current strategy $\mu_k$ and mean field $z_k,$ we would need to compute the stationary distribution of the model $P_{z_k,\mu_k},$ and set it as the next mean field.  This would preclude learning without a simulator, since running one step in the real world would immediately cause a mean field update to $z_{k+1},$ and induce a new model $P_{z_{k+1},.},$ making the system non-stationary.  

Since our RL algorithms only need a one-step McKean-Vlasov update under each strategy, an online learning approach at time $k$ when the underlying state distribution is $z_k,$ would be to apply $\mu_{k-1}$ to the system, with the resultant state distribution being $z_{k+1}.$  However, we face the issue that the samples obtained pertain to $P_{z_k,\mu_{k-1}},$ whereas the system model is now $P_{z_{k+1},.},$  and so an online sample-based trembling-hand strategy $\mu_{k}$ will lag the current system model by one step.  Fortunately, convergence of the TMFQ-learning approach is robust to this lag, and we present its online version in Algorithm~\ref{algo:ol-Q}.




\begin{algorithm} \label{algo:ol-Q}
	\caption{Online TMFQ-learning Algorithm for T-MFE}
	\begin{algorithmic}[1]
		\STATE Initialize mean field ${z}_0$ and strategy $\mu_{0}$ 
		\FOR {$k= 1, 2, \ldots $}
            \STATE Reset \textit{memory buffer}
		    \FOR {agents $i= 1, 2, \ldots, I$}
		        \STATE Given current state $s^{i}_{k},$ take action $a^{i}_{k} \sim \mu_{k-1}(s^{i}_{k}),$ get reward $r^{i}_{k},$ observe the next state
		        $s^{i}_{k+1}$
		        \STATE Add the sample $(s^{i}_{k}, a^{i}_{k}, r^{i}_{k}, s^{i}_{k+1})$ to the \textit{memory buffer}
		    \ENDFOR
            \STATE Perform TQ-learning on the \textit{memory buffer} to obtain $Q_{{z}_k}$ and ${\mu}_{k} = \pi^{\epsilon}_{Q_{z_k}}$
		\ENDFOR
	\end{algorithmic}
\end{algorithm}


Algorithm \ref{algo:ol-Q} aggregates the samples generated by executing strategy $\mu_{k-1}$ to estimate the TQ-value function, and then passes back $\mu_k$ to the agents.  We note that regeneration of a fraction of players and execution of trembling hand strategies ensures that we have sufficient samples of each state-action pair in the large population regime to ensure that off-policy learning such as TQ-learning on the \textit{memory buffer} converges to the optimal TQ-value function.   In order to characterize the sample complexity of Algorithm \ref{algo:ol-Q}, we employ bounds pertaining to synchronous Q-learning~\citep{even2003learning} using this guarantee that all state-action pairs are sampled a desired number of times.  The PAC result is similar to TMFQ-learning, with accuracy increasing in the number of players, rather than with the number of samples as in Theorem~\ref{thm:TQL-SC}.

\begin{theorem}\label{thm:sync-TQL-SC}
Let Assumption \ref{assump:Lips-R-P} hold. For any $0 \leq \bar{\epsilon}, \bar{\delta} < 1,$ let $k_{0} = k_{0}(\bar{\epsilon})$. In Algorithm \ref{algo:ol-Q}, for each $k \leq k_{0}$, assume that there are a total of $I = I_{0}$ number of players where $I_{0}$ is given as
\begin{align*}
\begin{split}
    I_{0} & = O \Bigg( \frac{|S| |A|}{\epsilon \zeta} \Bigg(\frac{B^{2} V^{2}_{\max}}{\beta^{2} \bar{\epsilon}^{2}} \ln \big( \frac{2 B k_{0}|\mathcal{S}||\mathcal{A}|V_{\max}}{\bar{\delta} \beta \bar{\epsilon}} \big)   \Bigg)^{\frac{1}{w}} \\
    & \quad + \Bigg( \frac{1}{\beta} \ln(\frac{B V_{\max}}{\bar{\epsilon}}) \Bigg)^{\frac{1}{1-w}} \Bigg),
\end{split}
\end{align*}
where $B = \frac{((1-\zeta) 2A)^{k_{0}} (C_{3}\epsilon)}{(1-\zeta) (A-1)}$, $A = \max \{ 1+C_{2}, C_{3}D \}$,  $D = (C_{1} + \gamma C_{2})/(1-\gamma)$, $V_{\max} = 1/(1-\gamma), \beta = (1-\gamma)/2$, $w \in (1/2, 1)$, $\zeta$ is a regeneration probability and $\epsilon$ is the trembling-hand strategy randomization.
Then,
\begin{align*}
    \mathbb{P}(\|z_{k_{0}} - z^{*} \| \leq 2 \bar{\epsilon}) \geq (1 - \bar{\delta}). 
\end{align*}
\end{theorem}
Note that a result similar to that of Corollary \ref{cor:cor-1} can easily be shown here as well under the same assumption. We omit details due to page limitations. 

\section{Experiments}
\label{sec:simulations}
We consider \emph{Infection Spread} described in Section \ref{sec:formulation}. State space is $\mathcal{S}=\{0,\ldots, 24\}$ with $0$ being the lowest health level. Action space is $\mathcal{A}=\{0, \ldots, 4\}$ with $4$ being the strongest preventive action, and $c_f$ is the infection intensity constant. Details of the parameters and more experiments are presented in the supplementary materials.  We evaluate the performance of our algorithms, TMFQ-Learning (TMFQ),  GMBL, and Online TMFQ-Learning (O-TMFQ) on this model, along with comparisons with Independent Q-Learning (IQL) and Mean Field Q-learning (MFQ) \citep{yang2018mean}.  In IQL,  each agent ignores other agents, maintains an individual Q-function and performs TQ-learning independently. We implement a variant of MFQ, where each agent maintains a Q-function parameterized by the average states of a subset of the population and performs TQ-Learning. We average over 20 runs in each experiment, and the dashed line and band in figures show the average and standard deviation, respectively.

Figure \ref{fig:CDF_Algos} shows the final mean field distribution obtained by each of our algorithms, simulated with $1000$ agents. Note that all of them converge to the same T-MFE, indicating the accuracy of O-TMFQ.  We next compare the performance O-TMFQ with  IQL and MFQ.  Figure \ref{fig:Algos_CDF} shows that the final equilibrium distributions of IQL and MFQ are inaccurate (not the true T-MFE), and that O-TMFQ results in higher states (health levels).  Figure \ref{fig:Algos_mean_state} shows the evolution of the mean health of the population, $\sum_{s} s z_{k}(s)$ with iteration number $k.$  The  mean health of the population quickly converges to a higher value under O-TMFQ while the corresponding value is lower under MFQ and IQL.


Figure \ref{fig:z_converge} shows rate of convergence of the mean field under different numbers of agents. Here, $z^*$ is the final mean field obtained by O-TMFQ.  We see that the asymptotically accurate mean field approximation becomes increasingly correct even with a relatively small number of agents of 500 or 1000.  We show the evolution of average health level of the population for O-TMFQ in Figure~\ref{fig:mean_state} for different values of $c_f$.  This plot indicates that the convergence of the algorithm is fairly fast.  Finally, Figure \ref{fig:CDF_OL} explores the impact of the mean field on the model and its associated equilibrium via $c_f$. As expected, more agents are in lower health states for larger $c_{f}$.

\section{Conclusions}

We introduced the notion of trembling-hand perfection to MFG as a means of providing strategically consistent exploration.  We showed existence of T-MFE in MFG with strategic complementarities, and developed an algorithm for computation.  Based on this algorithm, we developed model-free, model-based and fully online learning algorithms, and provided PAC bounds on their performance.  Experiments illustrated the accuracy and good convergence properties of our algorithms.

\clearpage
\bibliographystyle{abbrvnat}
\bibliography{ref}
\onecolumn

\clearpage
\appendix

\section{Algorithms Description}

\subsection{{Next-MF} Function for One-Step McKean-Vlasov Update}
\begin{algorithm} 
\label{algo:next-MF}
\caption{{Next-MF}}
	\begin{algorithmic}[1]
	    \STATE \textbf{Input:} Mean field $z$ and strategy $\mu$
		\STATE  Initialize $\hat{z}$, sample number  $j \leftarrow 0$
		\REPEAT
		    \STATE Sample the state $s \sim z(\cdot)$, action $a \in \mu(s, \cdot)$, next state $s' \sim P(\cdot | s, a, z)$
		    \STATE $\hat{z}(s') \leftarrow \hat{z}(s')  + 1$
		    \STATE $\bar{z}_{j} = \texttt{normalize}(\hat{z}$)
		    \STATE $j \leftarrow j+1$
		\UNTIL {$\norm{ \bar{z}_{j} -\bar{z}_{j-1}}_1 \leq \epsilon_2$}
		\STATE $\bar{z} = \bar{z}_{j}$
		\STATE \textbf{return $\bar{z}$}
	\end{algorithmic}
\end{algorithm}
Algorithm \ref{algo:next-MF} estimates the next mean field according to equation \eqref{eq:Mckean}. We maintain a frequency estimate of $s'$ in line 5 and normalize the frequency estimate to obtain a density function in line 6. Note that $s'$ is sampled according to line 4.

\section{Generative Model-Based Reinforcement Learning for T-MFE}
\label{sec:MBRL}
In this section we present a model-based variant of the T-BR algorithm.  Non-stationarity in our system (model changes at each step) implies that there is no single model that can be learned.  
Our approach follows the generation of a new model each time that the mean field evolves under a T-BR like approach.  

At each iteration, we first estimate the model $P(\cdot | \cdot, \cdot, z_{k})$ for the given $z_{k}$ using $n_{0}$ simulator samples  for each $(s,a)$.   We next define the approximate TQ-value operator $\widehat{F}_{z_{k}}$ as in \eqref{eq:F-operator} by replacing the actual model $P$ with the  estimated model $\widehat{P}$. It is straightforward to show that $\widehat{F}_{z_{k}}$ is also a contraction.  This ensures that approximate TQ-value iteration will converge to an approximate TQ-value function $\widehat{Q}^{*}_{z_{k}},$ which will, of course, have an error with respect to the true TQ-value function $Q^{*}_{z_{k}}$.  We determine the trembling-hand best response strategy with respect to  $\widehat{Q}^{*}_{z_{k}}$.  Finally, the next mean field $z_{k+1}$ is obtained using this strategy and the estimated model in the McKean-Vlasov update equation, denoted by $\hat{\Phi}(\cdot, \cdot)$. The GMBL algorithm is summarized in Algorithm~\ref{algo:MBL}.

Note that we are not estimating the model for all the possible mean fields, but only for the sequence $\{z_{k}\}$. So, if the process converges in a finite number of steps, then we need only a finite number of simulation samples. We show that this is indeed true in the theorem below. 

\begin{theorem}
\label{thm:MBL}
Let Assumption \ref{assump:Lips-R-P} hold. For any $0 \leq \bar{\epsilon}, \bar{\delta} < 1,$ let $k_{0} = k_{0}(\bar{\epsilon})$. In Algorithm \ref{algo:MBL}, for each $k \leq k_{0}$, assume that the estimate $\widehat{P}(\cdot | \cdot, \cdot, z_{k})$ is obtained by a total of $N_{0} = n_{0}|S||A|$ simulator samples where $n_{0}$ is given as
\begin{align*}
    n_{0} & = O \bigg( \max \bigg( \frac{2 V^{4}_{\max} B^{2}}{\bar{\epsilon}^{2}} \log \bigg( \frac{2 |\mathcal{S}| |\mathcal{A}| k_{0}}{\bar{\delta}} \bigg) , \nonumber \\ & \quad \quad \quad \frac{2 B^{2}}{\bar{\epsilon}^{2}} \log (\frac{2^{|\mathcal{S}|}  |\mathcal{S}| |\mathcal{A}| k_{0}}{\bar{\delta}} \bigg)\bigg)    \bigg)
\end{align*}
where $B = (1 + C_{2} + C_{3} D)^{k_{0}+1}$, $D = (C_{1} + \gamma C_{2})/(1-\gamma)$, $V_{\max} = 1/(1-\gamma)$. Then, 
\begin{align*}
    \mathbb{P}(\|z_{k_{0}} - z^{*} \| \leq 2 \bar{\epsilon}) \geq (1 - \bar{\delta}). 
\end{align*}
\end{theorem}

A result similar to that of Corollary \ref{cor:cor-1} can easily be shown here as well under the same assumption. The proof of theorem \ref{thm:MBL} is given in section \ref{appendix:learning-TMFE-offline} 

\begin{algorithm}[H]
\caption{GMBL Algorithm for T-MFE}
\label{algo:MBL}
\begin{algorithmic}[1]
\STATE Initialization: Initial mean field $z_{0}$ \\
\FOR {$k=0,1,2,..$}
\STATE For the mean field $z_{k}$ estimate the model to get $\widehat{P}(\cdot|\cdot,\cdot, z_{k})$ by taking\\ $n_{0}$ next-state samples for each state-action pair $(s,a)$
\STATE Compute  $\widehat{Q}^{*}_{z_{k}}$  using the approximate TQ-value iteration $\widehat{Q}_{m+1,z_{k}} = \widehat{F}_{z_{k}}(\widehat{Q}_{m,z_{k}})$
\STATE Compute the strategy 
${\mu}_{k} =  \pi^{\epsilon}_{\widehat{Q}^{*}_{z_{k}}}$. Then compute  the next mean field  $z_{k+1} = \widehat{\Phi}(z_{k}, \mu_{k}) $
\ENDFOR 
\end{algorithmic}
\end{algorithm}

\section{Proofs of Results in Section~\ref{sec:existence-comp}}
We use the following result from \cite{asadi2017alternative}
\begin{lemma}[\cite{asadi2017alternative}]
\label{lem:non-exp}
For any  $Q_{1}, Q_{2}$, and for any $s \in \mathcal{S},$ \[|G(Q_{1})(s) - G(Q_{2})(s)| \leq \max_{a} |Q_{1}(s, a) - Q_{2}(s, a) |. \]
\end{lemma}

\subsection{Proof of Proposition \ref{prop:F-mapping}}
\label{appendix:F-mapping}
\begin{proof}[Proof of Proposition \ref{prop:F-mapping}]
For any given  $(s,a) \in \mathcal{S} \times \mathcal{A}$ and mean field $z$,
\begin{align*}
    &|F_{z}(Q_{1})(s,a) -  F_{z}(Q_{2})(s,a)|  \le  \gamma |\sum_{s'} P(s' | s, a, z) (G(Q_1)(s') - G(Q_2)(s'))|  \\
    &\stackrel{(a)}{\leq}  \gamma \sum_{s'} P(s' | s, a, z)~ \max_{b} | Q_1(s', b) -  Q_2(s',b) | \leq \gamma \|Q_{1} - Q_{2} \|_{\infty}
\end{align*}
where $(a)$ follows from Lemma \ref{lem:non-exp}. Since $(s,a) \in \mathcal{S} \times \mathcal{A}$ was arbitrary, we have $\|F_{z}(Q_{1}) -  F_{z}(Q_{2})\|_{\infty} \leq \gamma  \|Q_{1} - Q_{2} \|_{\infty}$. Existence of a unique fixed point for $F_{z}$ follows directly from the  Banach's fixed point theorem since $F_{z}$ is a contraction. The claim that this unique fixed point is equal to $Q^{*}_{z}$ follows from the Bellman optimality principle. 
\end{proof}

\subsection{Proof of Theorem \ref{thm:tmfe-existence} and Theorem \ref{thm:T-MRD}}
\label{appendix:tmfe-existence}

We follow a proof approach that uses  the strategic complementarity conditions to establish some monotone properties, and exploit that to show the existence of T-MFE and the convergence of T-BR algorithm. We first state  some useful results from \cite{adlakha2013mean}.  Note that, however, since we are considering trembling-hand polices, proofs in  \cite{adlakha2013mean}  are not directly  applicable to our setting.

\begin{lemma}[Lemma 4 in \cite{adlakha2013mean}]
\label{lem:nondecreasing-expected-value}
Suppose that $V_{z}(s)$ is a non-decreasing bounded function in $s$ and has increasing differences in $s$ and $z$. Then, $\sum_{s' \in S} P(s'|s,a,z) V_{z}(s')$ is non-decreasing in $s$ and $a$ and has increasing differences in $(s,a)$ and $z$. Moreover, the function $V'_{z}(s)$ defined as, $V'_{z}(s) = \max_{a} ( r(s,a,z) + \gamma \sum_{s'} P(s'|s,a,z)V_{z}(s))$, is non-decreasing in $s$ and has  increasing differences in $s$ and $z$.
\end{lemma}
\begin{lemma}[Lemma 6 in \cite{adlakha2013mean}]
\label{lem:Omega-optimal-strategy}
Suppose that $V_{z}(s)$ is non-decreasing in s and has increasing differences in $s$ and $z$. Define a correspondence
\begin{align*}
\Omega(s,z) & = \arg \max\limits_{a \in \mathcal{A}} ~ ( r(s,a,z) + \gamma \sum_{s' \in S} P(s'|s,a,z) V_{z}(s') ).
\end{align*}
Then, $\Omega$ is a non-decreasing correspondence in $(s,z)$.
\end{lemma}

\begin{proposition}
\label{prop:opt-val-nonde-incre-diff}
Let $V^{*}_{z}$ and $\mu^{*}_{z}$ be the optimal trembling-hand value function and optimal trembling-hand strategy corresponding to mean field $z$. Then, $V^{*}_{z}(s)$ is non-decreasing in $s$ and has increasing differences in $s$ and $z$. Moreover, $\mu^{*}_{z}$  is stochastically non-decreasing in $s$ and $z$.
\end{proposition}
\begin{proof}
Let $\{Q_{m,z}, m \geq 0\}$ be the TQ-value iterates corresponding to a mean field $z$.  If $\arg \max_{b} Q_{m,z}(s, b)$ is not unique, then $\pi^{\epsilon}_{Q_{m,z}}$ can be defined in more than one  way. We define it as \begin{align}
\label{eq:eps-greedy-Q}
    \pi^{\epsilon}_{Q_{m, z}}(s,a) =
    \begin{cases}
        1 - \epsilon &\text{for } a = \sup \{\arg \max_{b} Q_{m,z}(s,b)\} \\
        \epsilon/(|\mathcal{A}|-1) &\text{for } a \not = \sup \{ \arg \max_{b} Q_{m,z}(s,b) \}
    \end{cases}
\end{align}
Note that the $\sup$ of a set is well defined with respect to a lattice according to the definition of MFG-SC. In the following we denote $ \pi^{\epsilon}_{Q_{m, z}}$ simply as $\mu_{m,z}$. 

By definition, $V^{*}_{z} = G(Q^{*}_{z})$ and define $V_{m, z} = G(Q_{m,z})$. To show that $V^{*}_{z}(s)$ is non-decreasing in $s$ and has increasing differences in $s$ and $z$, it suffices to show $V_{m,z}(s)$ has the same properties for all $m$. This is because, since $G$ is continuous, monotonicity and increasing differences are preserved under limits. 

Letting $Q_{0,z}(s,a) = 0$ for all $(s,a) \in \mathcal{S} \times \mathcal{A}$, we have  $Q_{1,z}(s,a) = r(s,a,z)$. Then, $Q_{1,z}(s,a)$ is non-decreasing in $s$ and has increasing differences in $s$ and $z$ by Definition \ref{def:mfg-sc}. Define the correspondence $\Omega_{m}$ as
\begin{align*}
    \Omega_{m}(s,z) & = \arg \max\limits_{a \in \mathcal{A}} ~ ( r(s,a,z) + \gamma \sum_{s' \in S} P(s'|s,a,z) V_{m-1,z}(s') ).
\end{align*}
By Lemma \ref{lem:Omega-optimal-strategy}, $\Omega_{1}(s,z) = \arg\max Q_{1,z}(s,a)$ is non-decreasing in $s$ and $z$. Using this, we can conclude that  $\mu_{1,z} = \pi^{\epsilon}_{Q_{1,z}}$  is stochastically non-decreasing in $s$ and $z$ for  $\epsilon$ such that $(1 - \epsilon) > \epsilon/(|\mathcal{A}|-1)$. To see this, first note that for $\epsilon = 0,$  the deterministic strategy $\pi^{\epsilon}_{Q_{1,z}}$ is non-decreasing in $s$ and $z$ by Lemma \ref{lem:Omega-optimal-strategy}. Now, for  $\epsilon$ with $(1 - \epsilon) > \epsilon/(|\mathcal{A}|-1),$ $\pi^{\epsilon}_{Q_{1,z}}$ is stochastically non-decreasing in $s$ and $z$ because the probability of the maximizing action is greater than all other actions. Thus, $\mu_{m,z}$ is stochastically non-decreasing in $s$ and $z$ for all $m \in \mathbb{N}$ if each $Q_{m,z}$ is non-decreasing in $s$ and has increasing differences in $s$ and $z$.

Now, recall that $V_{1,z}(s) = \sum_{a}\mu_{1,z}(s,a)Q_{1,z}(s,a) = (1 - \epsilon) Q_{1,z}(s,\bar{a}) + \epsilon \sum_{a}Q_{1,z}(s,a) $ where $\bar{a} = \sup \{\arg \max Q_{1,z}(s,a)\}$. Since $Q_{1,z}(s,a)$ is non-decreasing in $s$ and has increasing differences in $s$ and $z$ for all $a$, we can conclude that $V_{1,z}(s)$ is  is non-decreasing in $s$ and has increasing differences in $s$ and $z$.


As the induction hypothesis, suppose that both $Q_{m,z}(s,a)$ and $V_{m,z}(s)$ are non-decreasing in $s$ and $z$ and have increasing differences in $s$ and $z$ and $\mu_{m,z}$ is stochastically non-decreasing in $s$ and $z$.

Induction step is as follows: By Definition \ref{def:mfg-sc} and Lemma \ref{lem:nondecreasing-expected-value}, both $r(s,a,z)$ and $\gamma \sum_{s'} P(s'|s, a, z) V_{m,z}(s')$  are non-decreasing in $s$ and have increasing differences in $(s,a)$ and $z$. Therefore, $Q_{m+1,z}(s,a) = r(s,a,z) + \gamma \sum_{s'} P(s'|s, a, z) V_{m,z}(s')$ also satisfies the same properties. By Lemma \ref{lem:Omega-optimal-strategy}, notice that $\Omega_{m+1}(s,z) := \arg\max_{a \in \mathcal{A}} \{ Q_{m+1,z}(s,a) \}$ is non-decreasing in $(s,z)$ and therefore $\mu_{m+1,z} = \pi^{\epsilon}_{Q_{m+1,z}}$ defined as in \eqref{eq:eps-greedy-Q} is stochastically non-decreasing in $s$ and $z$ as argued before. Finally,  $V_{m+1,z}(s) = \sum_{a}\mu_{m+1,z}(s,a)Q_{m+1,z}(s,a) = (1 - \epsilon) Q_{m,z}(s,\bar{a}) + \epsilon \sum_{a}Q_{m,z}(s,a) $ where $\bar{a} =  \sup \{\arg \max Q_{m+1,z}(s,a)$\} is non-decreasing in $s$ and has increasing differences in $s$ and $z$ because $Q_{m+1,z}(s,a)$ is non-decreasing in $s$ and has increasing differences in $s$ and $z$ for all $a$.


Since TQ-value iteration converges, $Q^{*}_{z}(s,a)$ and $V^{*}_{z}(s)$ are non-decreasing in $s$ and have increasing differences in $s$ and $z$. By repeating the same argument again, $\Omega(s,z) = \arg\max_{a \in \mathcal{A}} \{ Q^{*}_{z}(s,a) \}$ is non-decreasing in $(s,z)$. So, $\mu^{*}_{z} = \pi^{\epsilon}_{Q^{*}_{z}}$ is stochastically non-decreasing in $s$ and $z$.
\end{proof}


Tarski’s fixed-point theorem \cite{tarski1955lattice} ensures that monotone functions on a lattice have a fixed point. We use that result to prove the existence of T-MFE. We first state Tarski’s fixed-point theorem for completeness. 

\begin{theorem}[Tarski's Fixed-point Theorem \cite{tarski1955lattice}]\label{thm:tarski}
Suppose that $\mathscr{L}$ is a nonempty complete lattice, and $T:\mathscr{L} \rightarrow \mathscr{L}$ is a non-decreasing function. Then the set of fixed points of $T$ is a nonempty complete lattice.
\end{theorem}

We now give the proof of Theorem \ref{thm:tmfe-existence}. 
\begin{proof}[{\bf Proof of Theorem \ref{thm:tmfe-existence}}]
For any strategy $\mu$ such that $\mu$ is stochastically non-decreasing in $s$, and for any given mean fields $z$ and $z'$, define the function $ K_{\mu,z}(z')$ as
\begin{align*}
    K_{\mu,z}(z')(s') = \sum_{s \in S} \sum_{a \in \mathcal{A}} z'(s) \mu(s,a) P(s'|s, a, z).
\end{align*}
From \cite[Lemma 7]{adlakha2013mean}, for $z_{2} \succeq_{SD} z_{1}$ and $z'_{2} \succeq_{SD} z'_{1}$ and $\mu_2(s,\cdot) \succeq_{SD} \mu_1(s,\cdot)$ with $\mu_2$ and $\mu_1$ are stochastically non-decreasing in $s$, we have $K_{\mu_{2},z_{2}}(z'_{2}) \succeq_{SD} K_{\mu_{1},z_{1}}(z'_{1})$. 

Observe that McKean-Vlasov update function $\Phi(z,\mu)$ is a special case of $K_{\mu,z}(z')$ by setting $z'=z$. By the above argument, $\Phi(z,\mu) = K_{\mu,z}(z)$ is stochastically non-decreasing in $(z,\mu)$ provided that $\mu$ is stochastically non-decreasing in $s$.   

Define $\Upsilon : \mathcal{P}(\mathcal{S}) \rightarrow \mathcal{P}(\mathcal{S})$ as $\Upsilon(z) = \Phi(z,\mu^{*}_{z})$ where $\mu^{*}_{z}$ is the optimal trembling-hand strategy corresponding to the mean field $z$. Recall that,  by Proposition \ref{prop:opt-val-nonde-incre-diff}, $\mu^{*}_{z}$ is stochastically non-decreasing in $s$ and $z$. Then,  for $z_{2} \succeq_{SD} z_{1}$, $\Upsilon(z_{2}) = \Phi(z_{2},\mu^{*}_{z_{2}}) \succeq_{SD} \Phi(z_{1},\mu^{*}_{z_{1}}) = \Upsilon(z_{1})$. From this, we can conclude that $\Upsilon$ is a stochastically non-decreasing function in $z$ and hence $\Upsilon(\cdot)$ has a fixed point by Tarski's theorem (Theorem \ref{thm:tarski}). In other words, there exists a $z^{*}$ such that $z^{*} = \Upsilon(z^{*}) = \Phi(z^{*}, \mu^{*}_{z^{*}})$. This implies that  there exists a mean field $z^{*}$ and strategy $\mu^{*}$ such that they satisfy the  optimality condition (i.e., $\mu^{*} = \mu^{*}_{z^{*}}$) and the consistency condition (i.e. $z^{*} = \Phi(z^{*}, \mu^{*}_{z^{*}})$). Thus, there exists a trembling-hand-perfect mean field equilibrium for mean field games with strategic complementarities.
\end{proof}

\begin{proof}[{\bf{Proof of Theorem \ref{thm:T-MRD}}}]
We exploit two key monotonicity properties established before. First, from proof of Theorem \ref{thm:tmfe-existence}, $\Upsilon(z)$ is stochastically non-decreasing in $z$. Second, from Proposition \ref{prop:opt-val-nonde-incre-diff}, an optimal trembling-hand strategy $\mu^{*}_{z}$ corresponding to a mean field $z$ is stochastically non-decreasing in $s$ and $z$.

Let $z_{0}$ be the smallest distribution by initialization in the $\preceq_{SD}$ ordering and let $\{Q_{z_{k}}, z_{k}, k \geq 0\}$ be the TQ-value functions and mean fields generated corresponding to Algorithm \ref{algo:T-MRD}. In the following, we denote an optimal trembling-hand strategy $\mu^{*}_{z_{k}}$ corresponding to a mean field $z_{k}$ simply as $\mu_{k}$. By Proposition \ref{prop:opt-val-nonde-incre-diff}, $\mu_{0} = \pi^{\epsilon}_{Q^{*}_{z_{0}}}$ is stochastically non-decreasing in $s$ and $z$. Hence we take the following as our induction base: $z_{0} \preceq_{SD} \Upsilon(z_{0}) = z_{1}$ and $\mu_{0}(s,\cdot) \preceq_{SD} \mu_{1}(s,\cdot)$ for all $s$ where $\mu_{1}$ is stochastically non-decreasing in $s$ and $z$.

Now as the induction hypothesis, suppose that $z_{0} \preceq_{SD} \Upsilon(z_{0}) = z_{1} \preceq_{SD} \Upsilon(z_{1}) = z_{2} \preceq_{SD} \cdots \preceq_{SD} \Upsilon(z_{k-1}) = z_{k}$ and that $\mu_{0}(s,\cdot) \preceq_{SD} \mu_{1}(s,\cdot) \preceq_{SD} \cdots \preceq_{SD} \mu_{k}(s,\cdot)$ for all $s$ where $\mu_{i}$ are stochastically non-decreasing in $s$ and $z$.

Then as an induction step, we have $\Upsilon(z_{k-1}) = z_{k} \preceq_{SD} z_{k+1} = \Upsilon(z_{k})$ because $\Upsilon(z)$ is stochastically non-decreasing in $z$. Now, since $ z_{k} \preceq_{SD} z_{k+1}$ and both $\mu_{k}$ and $\mu_{k+1}$ are stochastically non-decreasing in $s$ and $z$, it follows that $\mu_{k}(s,\cdot) \preceq_{SD} \mu_{k+1}(s,\cdot)$ for all $s$.

Observe that $\Pi^{\epsilon}$ is compact for a fixed $\epsilon > 0$ and since $(\mu_{k})_{k \in \mathbb{N}} \subset \Pi^{\epsilon}$ is a stochastically non-decreasing (monotone) sequence, there must be a pointwise limit $\mu^{*}$ such that $\mu_{k} \rightarrow \mu^{*}$ as $k \rightarrow \infty$. Moreover, $\mathcal{P}(\mathcal{S})$ is also compact since we assume that $|\mathcal{S}|$ is finite. Since $(z_{k})_{k \in \mathbb{N}} \subset \mathcal{P}(\mathcal{S})$ is a stochastically non-decreasing (monotone) sequence, there must be a limit $z^{*}$ such that $z_{k} \rightarrow z^{*}$ as $k \rightarrow \infty$.

It is straight forward to show that the optimal trembling-hand strategy $\mu^{*}_{z}$ and $\Upsilon(z)$ are continuous in $z$. So, since $\mu_{k} \rightarrow \mu^{*}_{z^{*}}$ and $z_{k} \rightarrow z^{*}$, we can conclude that $\mu^{*} = \mu^{*}_{z^{*}}$ and $z^{*} = \Upsilon(z^{*}) = \Phi(z^{*}, \mu^{*})$. 

This concludes the proof that T-BR converges to a T-MFE.
\end{proof}

\section{Proof of the Results in Section \ref{sec:learning-TMFE-offline}}
\label{appendix:learning-TMFE-offline}

\subsection{Proof of Theorem \ref{thm:T-MFE-convergence}}
\begin{proof}
We only sketch the proof since it is almost the same to the proof of Theorem \ref{thm:T-MRD}. At each time $k$ with a fixed $z_{k}$, the value of generalized Q-learning converges, i.e.,  $Q_{t,z_{k}} \rightarrow Q^*_{z_{k}}$ as $t \rightarrow \infty$. Then, under assumptions of the model, $Q^{*}_{z_{k}}$ satisfies the same complementarity properties for each $k \in \mathbb{N}$. Thus, we can conclude that, for all $k \in \mathbb{N}$, $\mu_{k}(s,\cdot) \preceq_{SD} \mu_{k+1}(s,\cdot)$ for all $s \in \mathcal{S}$ where each $\mu_{k}$ is stochastically non-decreasing in $s$ and $z$ and that $z_{k} \preceq_{SD} \Upsilon(z_k) = z_{k+1}$. The rest of the proof that the limits exist follows the same as in the proof of Theorem \ref{thm:T-MRD}.
\end{proof}

\subsection{Proof of Theorem \ref{thm:TQL-SC}}
We first prove some useful lemmas. 
\begin{lemma}
\label{lem:Q-ub-z}
Let $z_{1}$ and $z_{2}$ be two arbitrary mean fields and $Q^{*}_{z_{1}}$ and $Q^{*}_{z_{2}}$ be the optimal TQ-value functions corresponding to them. Then, under Assumption \ref{assump:Lips-R-P},
\begin{align}
    ||Q^{*}_{z_{1}}-Q^{*}_{{z}_{2}}||_{\infty} \le D ||z_{1} -z_{2}||_{1}
\end{align}
where $D = (C_{1} + \gamma C_{2})/(1- \gamma)^{2} $.
\end{lemma}
\begin{proof}
For any $(s,a) \in \mathcal{S} \times \mathcal{A}$,
\begin{align}
& |Q_{{z}_{1}}(s,a) - Q_{{z}_{2}}(s,a)|  \nonumber \\
&\le |r(s,a, z_{1}) - r(s,a, z_{2})| + \gamma |\sum_{s'}P(s'|s,a, z_{1}) G(Q_{z_{1}})(s') - \sum_{s'} P(s'|s,a, z_{2}) G(Q_{z_{2}})(s')| \nonumber \\
&\stackrel{(a)}{\leq} C_{1} \|z_{1} - z_{2} \|_{1} + \gamma |\sum_{s'}P(s'|s,a, z_{1}) |(G(Q_{z_{1}})(s') - G(Q_{z_{2}})(s')) | \nonumber  \\
&\hspace{3cm} + \gamma \sum_{s'} |P(s'|s,a, z_{1}) - P(s'|s,a, z_{2})|  |G(Q_{z_{2}})(s')| \nonumber \\
&\stackrel{(b)}{\leq} C_{1} \|z_{1} - z_{2} \|_{1} + \gamma |\sum_{s'}P(s'|s,a, z_{1}) ~ \max_{a} |Q_{z_{1}}(s', a) - Q_{z_{2}}(s', a) | \nonumber  \\
&\hspace{3cm} + \gamma  \sum_{s'} |P(s'|s,a, z_{1}) - P(s'|s,a, z_{2})|  |G(Q_{z_{2}})\|_{\infty} \nonumber \\
&\stackrel{(c)}{\leq} C_{1} \|z_{1} - z_{2} \|_{1} + \gamma \|Q_{z_{1}} - Q_{z_{2}} \|_{\infty}  + \frac{\gamma}{(1-\gamma)}  \|P(\cdot | \cdot, \cdot, z_{1}) - P(\cdot | \cdot, \cdot, z_{2}) \|_{1} \nonumber \\
\label{eq:Q-ub-z-1}
&\stackrel{(d)}{\leq} C_{1} \|z_{1} - z_{2} \|_{1} + \gamma \|Q_{z_{1}} - Q_{z_{2}} \|_{\infty}  + \frac{\gamma}{(1-\gamma)}  C_{2} \|z_{1} - z_{2} \|_{1}.
\end{align}
Here $(a)$ follows from Assumption 1.(i), $(b)$ from Lemma \ref{lem:non-exp}, $(c)$ from the fact that the maximum Q-value for a finite MDP is $1/(1-\gamma)$, and $(d)$ from Assumption 1.(ii). Since \eqref{eq:Q-ub-z-1}  is true for any $(s,a),$ by taking the maximum on the left hand side and re arranging, we get $
    \|Q_{z_{1}} - Q_{z_{2}} \|_{\infty}  \leq D \|z_{1} - z_{2} \|_{1},$
where $D = \frac{C_{1}}{1- \gamma} + \frac{\gamma C_{2}}{(1- \gamma)^{2}}$.
\end{proof}

\begin{lemma}
\label{lem:mf-evol-by-Q}
Let $Q_{1}, Q_{2}$ be two arbitrary TQ-value functions and let $\mu_{1}$ and $\mu_{2}$ be the trembling-hand strategies corresponding to them, i.e., $\mu_{1} = \pi^{\epsilon}_{Q_{1}}, \mu_{2} = \pi^{\epsilon}_{Q_{2}}$. Let $z_{1}, z_{2}$ be two arbitrary mean fields.  Then, under Assumption \ref{assump:Lips-R-P},
\begin{align}
   \| \Phi(z_{1}, \mu_{1}) -  \Phi(z_{2}, \mu_{2}) \|_{1} \leq (1+C_{2}) \|z_{1} - z_{2} \|_{1} + C_{3} \|Q_{1} - Q_{2} \|_{\infty}.
\end{align}
\end{lemma}
\begin{proof}
We have,
\begin{align*}
     &\| \Phi(z_{1}, \mu_{1}) -  \Phi(z_{2}, \mu_{2}) \|_{1} = \sum_{s'}|\Phi(z_{1}, \mu_{1})(s') -  \Phi(z_{2}, \mu_{2})(s')|  \\
    &= \sum_{s'} | \sum_{s \in \mathcal{S}} \sum_{a \in \mathcal{A}} z_{1}(s) \mu_{1}(s, a) P(s' |s, a, z_{1}) - \sum_{s \in \mathcal{S}} \sum_{a \in \mathcal{A}} z_{2}(s) \mu_{2}(s, a) P(s' |s, a, z_{2}) | \\
    &= \sum_{s'}  |\sum_{s \in \mathcal{S}} z_{1}(s) P_{z_{1}, \mu_{1}}(s'|s) - \sum_{s \in \mathcal{S}} z_{2}(s) P_{z_{2}, \mu_{2}}(s'|s) | \\
    &\leq \sum_{s'}  |\sum_{s \in \mathcal{S}} (z_{1}(s) - z_{2}(s)) P_{z_{1}, \mu_{1}}(s'|s) | + \sum_{s'} |  \sum_{s \in \mathcal{S}} z_{2}(s) ( P_{z_{1}, \mu_{1}}(s'|s) - P_{z_{2}, \mu_{2}}(s'|s) ) | \\
    &\leq \|z_{1} - z_{2} \|_{1} + \|P_{z_{1}, \mu_{1}} - P_{z_{2}, \mu_{2}} \|_{1} \\
    &\leq \|z_{1} - z_{2} \|_{1} + \|P_{z_{1}, \mu_{1}} - P_{z_{2}, \mu_{1}} \|_{1} + \|P_{z_{2}, \mu_{1}} - P_{z_{2}, \mu_{2}} \|_{1}  \\
    &\stackrel{(a)}{\leq} \|z_{1} - z_{2} \|_{1} + C_{2}  \|z_{1} - z_{2} \|_{1} + C_{3} \|Q_{1} -Q_{2} \|_{\infty}
\end{align*}
where $(a)$ follows from Assumption \ref{assump:Lips-R-P}. 
\end{proof}

We use the following Q-learning sample complexity results from \cite{even2003learning}.
\begin{theorem}[Theorem 4 in \cite{even2003learning}]
\label{thm:QL-SC}
Let $Q_{t}$ be the $t$-th update in Q-learning algorithm using a polynomial time learning rate given as $\alpha_{t}(s, a) = {1}/({n_{t}(s,a)}+1)^{w}$, where $n_{t}(s,a)$ is the number of times the state-action pair $(s,a)$ is visited until time $t$ and $w \in (1/2, 1)$. Let $L$ be the upper bound on the covering time. Then, $ \mathbb{P}(\norm{Q_{T_{0}} - Q^{*}}_{\infty} \le \epsilon_{3}) \geq (1 - \delta_{3})$, for any $0 < \epsilon_{3}, \delta_{3} < 1,$ given that
\begin{align}
\label{eq:asynch-Q-sample}
    T_{0} = O \Bigg( \Bigg( \frac{ L^{1+3w} V^{2}_{\max} \ln \big( \frac{|\mathcal{S}||\mathcal{A}|V_{\max}}{\delta_{3} \beta \epsilon_{3}} \big)}{\beta^{2} \epsilon^{2}_{3}} \Bigg)^{\frac{1}{w}} + \Bigg( \frac{L}{\beta} \ln(\frac{V_{\max}}{\epsilon_{3}}) \Bigg)^{\frac{1}{1-w}} \Bigg).
\end{align}
where $V_{\max} = 1/(1-\gamma), \beta = (1-\gamma)/2$.
\end{theorem}
Here the covering time of a state-action pair sequence is  the number of steps needed to visit all state-action pairs starting from any arbitrary state-action pair. 

We note that the Q-learning update used in TQ-learning algorithm satisfies all the conditions necessary for the above theorem. So, we will use the above result. We refer the reader to \cite{even2003learning} for the details.

We now give the proof Theorem \ref{thm:TQL-SC}.
\begin{proof}[{\bf Proof of Theorem \ref{thm:TQL-SC}} ]
Let $\{z_{k} \}$ and  $\{ \mu_{k}\}$ be the sequences of  mean fields and strategies generated by the TMFQ-learning algorithm. Let $\{Q_{t,k}, t \geq 0\}$ be TQ-learning iterates corresponding to the mean field $z_{k}$ and let $Q_{k} = Q_{T_{0}, k}$ where $T_{0}$ is as given in \eqref{eq:asynch-Q-sample}. We assume that the number of samples used in the \texttt{Next-MF} function is such that, for any given mean field $z$ and strategy $\mu$, \texttt{Next-MF} function returns a mean field $z'$ such that $\mathbb{P}(\|z' - \Phi(z, \mu)\|_{1} \leq \epsilon_{3}) \geq (1 - \delta_{3}).$

Define the event  $E_{k} = \{ \|Q_{k} - Q^{*}_{z_{k}} \|_{\infty} \leq \epsilon_{3}~\text{and}~ \|z_{k+1} - \Phi(z_{k}, \mu_{k})\|_{1} \leq \epsilon_{3} \}$. Then, according to Theorem \ref{thm:QL-SC} and the assumption on the \texttt{Next-MF} function, $\mathbb{P}(E_{k}) \geq (1 - 2 \delta_{3})$. Define the event $E = \cap^{k_{0}}_{k=1} E_{k}$. So, $\mathbb{P}(E) \geq (1 - 2 k_{0} \delta_{3})$.  We will now  analyze the TMFQ-learning algorithm conditioned on the event $E$. 

Let $\{\bar{z}_{k} \}$ and  $\{\bar{\mu}_{k}\}$ be the sequences of  mean fields and strategies generated by the T-BR  algorithm. We assume that T-BR algorithm and TMFQ-learning algorithm have the same initialization, i.e.,  $\bar{z}_{0} = z_{0}$. Now,  conditioned on the event $E$,
\begin{align*}
    \|\bar{z}_{k+1} - z_{k+1} \|_{1} & \leq \| \Phi(\bar{z}_{k}, \bar{\mu}_{k}) -  \Phi({z}_{k}, {\mu}_{k}) \|_{1} +  \|  \Phi({z}_{k}, {\mu}_{k}) - z_{k+1} \|_{1} \\
    & \stackrel{(a)} \leq (1+C_{2}) \|\bar{z}_{k} - {z}_{k} \|_{1} + C_{3} \|Q^{*}_{\bar{z}_{k}} - Q^{*}_{{k}} \|_{\infty} +  \epsilon_{3} \\
    & \leq (1+C_{2}) \|\bar{z}_{k} - {z}_{k} \|_{1} + C_{3} \|Q^{*}_{\bar{z}_{k}} - Q^{*}_{z_{k}} \|_{\infty} + C_{3} \|Q^{*}_{z_{k}} - Q^{*}_{k}  \|_{\infty}+ \epsilon_{3} \\
    & \stackrel{(b)} \leq (1+C_{2}) \|\bar{z}_{k} - {z}_{k} \|_{1} + C_{3} D   \|\bar{z}_{k} - z_{k} \|_{1} + (C_{3} +1) \epsilon_{3}  \\
    & \leq (1 + C_{2} + C_{3} D) \|\bar{z}_{k} - z_{k} \|_{1} + (C_{3} +1) \epsilon_{3}
\end{align*}
Here $(a)$ follows from Lemma \ref{lem:mf-evol-by-Q} and the assumption on the \texttt{Next-MF} function and $(b)$ follows from  Lemma \ref{lem:Q-ub-z}. 

Iteratively applying the above inequality, we get  $ \|\bar{z}_{k_{0}} - z_{k_{0}} \|_{1} \leq B \epsilon_{3},$ where $B = (1 + C_{2} + C_{3} D)^{k_{0}+1} (C_{3}+1)$.  Now, $\|\bar{z}^{*} - z_{k_{0}} \|_{1} \leq \|\bar{z}^{*} - \bar{z}_{k_{0}} \|_{1} + \|\bar{z}_{k_{0}} - z_{k_{0}} \|_{1} \leq \bar{\epsilon} + B \epsilon_{3}$ because $\|\bar{z}^{*} - \bar{z}_{k_{0}} \|_{1} \leq \bar{\epsilon}$ by the definition of $k_{0}$. 

So, $\mathbb{P}( \|\bar{z}^{*} - z_{k_{0}} \|_{1} \leq \bar{\epsilon} + B \epsilon_{3}) \geq \mathbb{P}(E) = 1 - 2 k_{0} \delta$.

Setting $ \epsilon_{3} = \bar{\epsilon}/B$ and $\delta_{3} = \delta/2k_{0}$, and using the corresponding $T_{0}$ from \eqref{eq:asynch-Q-sample}, we get the desired result. 
\end{proof}

\begin{proof}[{\bf Proof of Corollary \ref{cor:cor-1}} ]
Conditioned on the event $E$ as defined in the proof Theorem \ref{thm:TQL-SC}, we get, 
\begin{align*}
    \|\bar{z}_{k+1} - z_{k+1} \|_{1} & \leq \| \Phi(\bar{z}_{k}, \bar{\mu}_{k}) -  \Phi({z}_{k}, {\mu}_{k}) \|_{1} +  \|  \Phi({z}_{k}, {\mu}_{k}) - z_{k+1} \|_{1} \\
    &\stackrel{(a)} \leq C_{4} \|\bar{z}_{k} - {z}_{k} \|_{1} + C_{5} \|Q^{*}_{\bar{z}_{k}} - Q^{*}_{{k}} \|_{\infty} +  \epsilon_{3} \\
     & \leq C_{4} \|\bar{z}_{k} - {z}_{k} \|_{1} + C_{5} \|Q^{*}_{\bar{z}_{k}} - Q^{*}_{z_{k}} \|_{\infty} + C_{5} \|Q^{*}_{z_{k}} - Q^{*}_{k}  \|_{\infty}+ \epsilon_{3} \\
    &\stackrel{(b)} \leq C_{4} \|\bar{z}_{k} - {z}_{k} \|_{1} + C_{5} D   \|\bar{z}_{k} - z_{k} \|_{1} + (C_{5} +1) \epsilon_{3}  \\
    & \leq (C_{4} + C_{5} D) \|\bar{z}_{k} - z_{k} \|_{1} + (C_{5} +1) \epsilon_{3}
\end{align*}
Here $(a)$ follows from Assumption \ref{assump:MV-strong} and  and $(b)$ follows from  Lemma \ref{lem:Q-ub-z}. 

Iteratively applying the above inequality, we get  $ \|\bar{z}_{k_{0}} - z_{k_{0}} \|_{1} \leq B \epsilon_{3},$ where\\ $B = (C_{5}+1)/(1- ((C_{4} + C_{5} D))) $.  Rest of the proof is similar to that of Theorem \ref{thm:TQL-SC}.
\end{proof}

\subsection{Proof of Theorem \ref{thm:MBL}}
For a given mean field $z,$ let $\widehat{P}(\cdot|\cdot ,\cdot , z)$ be the estimate of the model obtained by taking $n_{0}$ next-state samples  for each $(s, a)$.  Let $\widehat{F}_{z}$ be the approximate TQ-value operator obtained by replacing $P$ by $\widehat{P}$ in \eqref{eq:F-operator}. Similar to the proof of Proposition \ref{prop:F-mapping}, it is straight forward to show that $\widehat{F}_{z}$ is a contraction. Let $\widehat{Q}^{*}_{z}$ be its unique fixed point. Since $\widehat{P}$ is different from $P$,  $\widehat{Q}^{*}_{z}$ and ${Q}^{*}_{z}$  will also be different. However, for sufficiently large $n_{0},$ we can give the following bound.  
\begin{lemma}
\label{lem:model-L1}
For any $0 < \epsilon_{4}, \delta_{4} < 1,$ 
\begin{align}
    \mathbb{P}(\| \widehat{Q}^{*}_{z} - {Q}^{*}_{z} \|_{\infty} \leq \epsilon_{4})  \geq ( 1- \delta_{4}),~~\text{for}~~  n_{0} \geq \frac{2 V^{4}_{\max}}{\epsilon^{2}_{4}} \log \left(\frac{2 |\mathcal{S}| |\mathcal{A}|}{\delta_{4}} \right)
\end{align}
where $V_{\max} = 1/(1-\gamma)$. 
\end{lemma}
\begin{proof}
\begin{align}
    &| \widehat{Q}^{*}_{z}(s, a) - Q^{*}_{z}(s, a) | =  | \widehat{F}_{z}(\widehat{Q}^{*}_{z})(s, a) - F_{z}(Q^{*}_{z})(s, a) | \nonumber \\
    &=  \gamma |\sum_{s'}  \widehat{P}(s'|s,a, z) G(\widehat{Q}^{*}_{z})(s') - \sum_{s'} P(s'|s,a, z) G(Q^{*}_{z})(s')| \nonumber \\
     & \leq  \gamma |\sum_{s'}  \widehat{P}(s'|s,a, z) (G(\widehat{Q}^{*}_{z})(s') -  G(Q^{*}_{z})(s') ) | +  \gamma | \sum_{s'}  (\widehat{P}(s'|s,a, z) - P(s'|s,a, z) )  G(Q^{*}_{z})(s')| \nonumber \\
     \label{eq:mbl-step1}
      &\stackrel{(a)}{\leq}  \gamma \|\widehat{Q}^{*} - Q^{*}_{z} \|_{\infty} + \gamma | \sum_{s'}  (\widehat{P}(s'|s,a, z)) - P(s'|s,a, z) )  G(Q^{*}_{z}(s')|,
\end{align}
where $(a)$ follows from Lemma \ref{lem:non-exp}. 

For bounding $|\sum_{s'}  (\widehat{P}(s'|s,a, z) - P(s'|s,a, z) )  G(Q^{*}_{z})(s')|$, note that $\sum_{s'}  \widehat{P}(s'|s,a, z)) G(Q^{*}_{z})(s')$ is  an unbiased estimated of $\sum_{s'}  {P}(s'|s,a, z) G(Q^{*}_{z})(s')$. Also note that \\ $\max_{s'} |G(Q^{*}_{z})(s')| \leq  1/(1-\gamma) = V_{\max}$. So, by applying Hoeffding's inequality, for a given $(s,a),$ we get
\begin{align*}
    \mathbb{P}(|\sum_{s'}  (\widehat{P}(s'|s,a, z)) - P(s'|s,a, z) )  G(Q^{*}_{z})(s') | \geq \epsilon ) \leq 2 ~\text{exp}\left( \frac{-n_{0} \epsilon^{2}}{2 V^{2}_{\max}} \right).
\end{align*}
Using the union bound argument, for all $(s, a),$ we get
\begin{align*}
    \mathbb{P}(|\sum_{s'}  (\widehat{P}(s'|s,a, z)) - P(s'|s,a, z) )  G(Q^{*}_{z})(s') | \leq \epsilon ) \geq 1 - 2 |\mathcal{S}| |\mathcal{A}| ~\text{exp}\left( \frac{-n_{0} \epsilon^{2}}{2 V^{2}_{\max}} \right). 
\end{align*}
So, with  $n_{0} \geq \frac{2 V^{2}_{\max}}{\epsilon^{2}} \log \left(\frac{2 |\mathcal{S}| |\mathcal{A}|}{\delta} \right)$,
\begin{align}
\label{eq:ub-H-1}
    \mathbb{P}(|\sum_{s'}  (\widehat{P}(s'|s,a, z)) - P(s'|s,a, z) )  G(Q^{*}_{z})(s') | \leq \epsilon )  \leq 1 -\delta,~\forall{(s, a) \in \mathcal{S} \times \mathcal{A}}. 
\end{align}
Now, for the above $n_{0}$, from \eqref{eq:mbl-step1} and \eqref{eq:ub-H-1},  $ | \widehat{Q}^{*}_{z}(s, a) - Q^{*}_{z}(s, a)| \leq \gamma \|\widehat{Q}^{*} - Q^{*}_{z} \|_{\infty} + \epsilon,$ with a probability greater than $(1-\delta)$ for all $(s,a)$. This implies that $\|\widehat{Q}^{*} - Q^{*}_{z} \|_{\infty} \leq \epsilon/(1-\gamma) = \epsilon V_{\max}$. Now, using $\epsilon = \epsilon_{4}/V_{\max}$ and $\delta = \delta_{4}$ in the expression for $n_{0}$ above, we get the desired result. 
\end{proof}

We now bound the error in the Mckean-Vlasov update due to of replacing $P$ by $\widehat{P}$. 
\begin{lemma}
\label{lem:McKean-Hat}
Let $\widehat{\Phi}$ be the approximate  Mckean-Vlasov update function as defined in \eqref{eq:Mckean} but  by replacing $P$ by $\widehat{P}$. Then, 
\begin{align}
\label{eq:Phi-hat-bd-1}
\mathbb{P}(\|\widehat{\Phi}(z, \mu) - {\Phi}(z, \mu)\|_{1} \leq \epsilon_{4}) \geq (1 - \delta_{4}),~ \text{for}~ n_{0} \geq \frac{2}{\epsilon^{2}_{4}} \log\left(\frac{2^{|\mathcal{S}|}  |\mathcal{S}| |\mathcal{A}|}{\delta_{4}} \right)
\end{align}
\end{lemma}
\begin{proof}
From \cite{Weissman03}, for a given $(s,a)$, 
\begin{align*}
    \mathbb{P}(\|\widehat{P}(\cdot | s, a, z) - {P}(\cdot | s, a, z)\|_{1} \geq \epsilon_{4}) \leq 2^{|\mathcal{S}|} \text{exp}\left(\frac{-n \epsilon^{2}_{4}}{2} \right).
\end{align*}
By the union bound argument, for all $(s,a)$, we get
\begin{align*}
    \mathbb{P}(\|\widehat{P}(\cdot | s, a, z) - {P}(\cdot | s, a, z)\|_{1} \leq \epsilon_{4}) \geq 1 - |\mathcal{S}| |\mathcal{A}| 2^{|\mathcal{S}|} \text{exp}\left(\frac{-n \epsilon^{2}_{4}}{2} \right).
\end{align*}
So, with $n_{0} \geq \frac{2}{\epsilon^{2}_{4}} \log\left(\frac{2^{|\mathcal{S}|}  |\mathcal{S}| |\mathcal{A}|}{\delta_{4}} \right)$
\begin{align}
\label{eq:Phi-hat-bd-2}
    \mathbb{P}(\|\widehat{P}(\cdot | s, a, z) - {P}(\cdot | s, a, z)\|_{1} \leq \epsilon_{4}) \geq 1 - \delta_{4}, ~\forall{(s, a) \in \mathcal{S} \times \mathcal{A}}. 
\end{align}

Now, with the above $n_{0}$, with a probability greater than $(1 - \delta_{4})$, we get
\begin{align*}
    \|\widehat{\Phi}(z, \mu) - {\Phi}(z, \mu)\|_{1} &= \sum_{s'} |\widehat{\Phi}(z, \mu)(s') - {\Phi}(z, \mu)(s')|\\
    &\leq  \sum_{s'} \sum_{s} \sum_{a} z(s) \mu(s,a) | \widehat{P}(s' | s, a, z) - {P}(s' | s, a, z) | \\
    &= \sum_{s} \sum_{a} z(s) \mu(s,a) \| \widehat{P}(\cdot | s, a, z) - {P}(\cdot | s, a, z) \|_{1} \leq  \epsilon_{4}
\end{align*}
where the last inequality follows from \eqref{eq:Phi-hat-bd-2} and the fact that $\sum_{s} \sum_{a} z(s) \mu(s,a) = 1$. 
\end{proof}

We now give the proof Theorem \ref{thm:MBL}.
\begin{proof}[{\bf Proof of Theorem \ref{thm:MBL}}]

Let $\{z_{k} \}$ and  $\{ \mu_{k}\}$ be the sequences of  mean fields and strategies generated by the GMBL algorithm. Let $n_{0}$ be the maximum of the two values given by Lemma \ref{lem:model-L1} and Lemma \ref{lem:McKean-Hat}, i.e,  
\begin{align*}
    n_{0} = \max \left( \frac{2 V^{4}_{\max}}{\epsilon^{2}_{4}} \log \left(\frac{2 |\mathcal{S}| |\mathcal{A}|}{\delta_{4}} \right) , \frac{2}{\epsilon^{2}_{4}} \log\left(\frac{2^{|\mathcal{S}|}  |\mathcal{S}| |\mathcal{A}|}{\delta_{4}} \right)\right) 
\end{align*}

Define the event  $E_{k} = \{\| \widehat{Q}^{*}_{z_{k}} - {Q}^{*}_{z_{k}} \|_{\infty} \leq \epsilon_{4} ~\text{and}~  \|\widehat{\Phi}(z_{k}, \mu_{k}) - {\Phi}(z_{k}, \mu_{k})\|_{1} \leq \epsilon_{4}\}$. Then, from Lemma \ref{lem:model-L1} and Lemma \ref{lem:McKean-Hat}, we get  $\mathbb{P}(E_{k}) \geq (1 - 2 \delta_{4})$. Define the event $E = \cap^{k_{0}}_{k=1} E_{k}$. So, $\mathbb{P}(E) \geq (1 - 2 k_{0} \delta_{4})$. We will now analyze the GMBL algorithm conditioned on the event $E$. 

Let $\{\bar{z}_{k} \}$ and  $\{\bar{\mu}_{k}\}$ be the sequences of  mean fields and strategies generated by the T-BR  algorithm. We assume that T-BR algorithm and GMBL algorithm have the same initialization, i.e.,  $\bar{z}_{0} = z_{0}$. Now, conditioned on the event $E$
\begin{align*}
    \|\bar{z}_{k+1} - z_{k+1} \|_{1} &=   \| \Phi(\bar{z}_{k}, \bar{\mu}_{k}) -  \widehat{\Phi}({z}_{k}, {\mu}_{k}) \|_{1} \\
    &\leq  \| \Phi(\bar{z}_{k}, \bar{\mu}_{k}) -  {\Phi}({z}_{k}, {\mu}_{k}) \|_{1}  + \| {\Phi}({z}_{k}, {\mu}_{k}) -  \widehat{\Phi}({z}_{k}, {\mu}_{k}) \|_{1} \\
    &\stackrel{(a)}{\leq} (1+C_{2}) \|\bar{z}_{k} - {z}_{k} \|_{1} + C_{3} \|Q^{*}_{\bar{z}_{k}} - \widehat{Q}^{*}_{z_{k}} \|_{\infty} + \epsilon_{4} \\
    &\leq (1+C_{2}) \|\bar{z}_{k} - {z}_{k} \|_{1} + C_{3} \|Q^{*}_{\bar{z}_{k}} - {Q}^{*}_{z_{k}} \|_{\infty} + C_{3} \|{Q}^{*}_{z_{k}} - \widehat{Q}^{*}_{z_{k}} \|_{\infty} +\epsilon_{4} \\
    & \stackrel{(b)}{\leq} (1+C_{2}) \|\bar{z}_{k} - {z}_{k} \|_{1} + C_{3} D   \|\bar{z}_{k} - z_{k} \|_{1} + (C_{3} +1) \epsilon_{4}  \\
    &\leq (1 + C_{2} + C_{3} D) \|\bar{z}_{k} - z_{k} \|_{1} + (C_{3} +1)  \epsilon_{4}
\end{align*}
Here $(a)$ follows from Lemma \ref{lem:mf-evol-by-Q} and $(b)$ follows from Lemma \ref{lem:Q-ub-z}. 

Iteratively applying the above inequality, we get  $ \|\bar{z}_{k_{0}} - z_{k_{0}} \|_{1} \leq B \epsilon_{4},$ where $B = (1 + C_{2} + C_{3} D_{1})^{k_{0}+1} (C_{3}+1)$.  Now, conditioned on the event $E$, we have
\begin{align*}
     \|\bar{z}^{*} - z_{k_{0}} \|_{1} \leq \|\bar{z}^{*} - \bar{z}_{k_{0}} \|_{1} + \|\bar{z}_{k_{0}} - z_{k_{0}} \|_{1} \leq \bar{\epsilon} + B \epsilon_{4}
\end{align*}
because $\|\bar{z}^{*} - \bar{z}_{k_{0}} \|_{1} \leq \bar{\epsilon}$ by the definition of $k_{0}$.

Setting $ \epsilon_{4} = \bar{\epsilon}/B$ and $\delta_{4} = \delta/2k_{0}$ in the expression for $n_{0},$ we get the desired result. 
\end{proof}

\section{Proof of the Results in Section \ref{sec:learning-TMFE-online}}
\label{appendix:learning-TMFE-onine}

\subsection{Proof of Theorem \ref{thm:sync-TQL-SC}}
While Q-learning algorithm is an asynchronous process since a particular state-action pair is updated at a time, if all state-action pairs are updated at each time, it is called synchronous Q-learning algorithm \cite{even2003learning} and we define synchronous TQ-learning algorithm as follows: For a fixed mean field $z$,
\begin{align}\label{eq:sync-TQ-update}
    & Q_{0,z}(s,a) = 0 ~~\text{for all}~~ (s,a) \in S \times A \nonumber \\
    & Q_{t+1,z}(s,a) = (1-\alpha)Q_{t,z}(s,a) + \alpha (r(s,a,z) + \gamma G(Q)(s)) ~~ \forall~ (s,a) \in S \times A
\end{align}
where $\alpha_{t}$ is the appropriate learning rate. It can be shown that $Q_{t+1,z} \rightarrow Q^{*}_{z}$ as $t \rightarrow \infty$ and observe that the synchronous TQ-value function update \eqref{eq:sync-TQ-update} requires that all state-action pairs are sampled at each $t$. Due to the availability of a large population of players that regenerate to occupy all states, and explore all actions via trembling hand strategies, this is a mild condition.

We use the following synchronous Q-learning sample complexity results from \cite{even2003learning}.
\begin{theorem}[Theorem 2 in \cite{even2003learning}]
\label{thm:sync-QL-SC}
Let $Q_{t}$ be the $t$-th update in synchronous Q-learning algorithm using a polynomial time learning rate given as $\alpha_{t}(s, a) = {1}/({n_{t}(s,a)}+1)^{w}$, where $n_{t}(s,a)$ is the number of times the state-action pair $(s,a)$ is visited until time $t$ and $w \in (1/2, 1)$. Then, $ \mathbb{P}(\norm{Q_{I_{0}} - Q^{*}}_{\infty} \le \epsilon_{3}) \geq (1 - \delta_{3})$, for any $0 < \epsilon_{3}, \delta_{3} < 1,$ given that
\begin{align}
\label{eq:synch-Q-sample-1}
    I_{0} = O \Bigg( \Bigg( \frac{ V^{2}_{\max} \ln \big( \frac{|\mathcal{S}||\mathcal{A}|V_{\max}}{\delta_{3} \beta \epsilon_{3}} \big)}{\beta^{2} \epsilon^{2}_{3}} \Bigg)^{\frac{1}{w}} + \Bigg( \frac{1}{\beta} \ln(\frac{V_{\max}}{\epsilon_{3}}) \Bigg)^{\frac{1}{1-w}} \Bigg).
\end{align}
where $V_{\max} = 1/(1-\gamma), \beta = (1-\gamma)/2$.
\end{theorem}

Unlike asynchronous Q-learning sample complexity, there is no dependence of covering time since all state-action pairs are updated at each time. We note that the synchronous Q-learning update used in synchronous TQ-learning algorithm satisfies all the conditions necessary for the above theorem. So, we will use the above result. We refer the reader to \cite{even2003learning} for the details.

We provide a slight modification of T-BR algorithm which intentionally include a mismatch observed in online TMFQ-learning algorithm. As done in the proof of Theorem \ref{thm:TQL-SC} where we compare trajectories of T-BR and TMFQ-learning algorithms, we compare modified T-BR and Online TMFQ-learning in a similar way because there is no loss of generality, i.e., modified T-BR also converges to a T-MFE. We can compare with T-BR and Online TMFQ-learning but there is an additional term. To elaborate about the Modified T-BR Algorithm \ref{algo:modified-T-MRD}, at each $k$ with $z_{k}$, Algorithm \ref{algo:modified-T-MRD} computes $Q^{*}_{z_{k}}$ and $\mu_{k}$ but, in McKean-Vlasov equation, $\mu_{k-1}$ is employed that is computed in previous step $k-1$, i.e., $\Phi(z_{k}, \mu_{k-1})$, and uses $\mu_{k}$ in next time step $k+1$. This is intentional and almost sure convergence to a T-MFE can be proved similarly.
\begin{algorithm}
\caption{Modified T-BR Algorithm}
\label{algo:modified-T-MRD}
\begin{algorithmic}[1]
\STATE Initialization: Initial mean field $z_{0}$ and strategy $\mu_{0}$ 
\FOR {$k=1,2,3, \ldots $}
\STATE For the mean field $z_{k}$, compute the optimal TQ-Value function $Q^{*}_{z_{k}}$ using the TQ-value iteration $Q_{m+1,z_{k}} = F(Q_{m,z_{k}})$
\STATE Compute the strategy $\mu_{k} = \Psi(z_{k})$ as the trembling-hand strategy w.r.t $Q^{*}_{z_{k}}$, i.e, 
$\mu_{k} =  \pi^{\epsilon}_{Q^{*}_{z_{k}}}$
\STATE Compute the next mean field  $z_{k+1} = \Phi(z_{k}, \mu_{k-1})$
\ENDFOR 
\end{algorithmic}
\end{algorithm}


\begin{proof}[{\bf Proof of Theorem \ref{thm:sync-TQL-SC}}]
Let $\{z_{k} \}$ and  $\{ \mu_{k}\}$ be the sequences of mean fields and strategies generated by the Online TMFQ-learning algorithm. Let $\{Q_{t,k}, t \geq 0\}$ be Offline (or Batch) TQ-learning iterates corresponding to the mean field $z_{k}$ and let $Q_{k} = Q_{I_{0}, k}$ where $I_{0}$ is as given in \eqref{eq:synch-Q-sample-1}.

Define the event $E_{k} = \{ \|Q_{k} - Q^{*}_{z_{k}} \|_{\infty} \leq \epsilon_{3} \}$. Then, according to Theorem \ref{thm:sync-QL-SC}, $\mathbb{P}(E_{k}) \geq (1 - 2 \delta_{3})$. Define the event $E = \cap^{k_{0}}_{k=1} E_{k}$. So, $\mathbb{P}(E) \geq (1 - 2 k_{0} \delta_{3})$.  We will now analyze the Online TMFQ-learning algorithm conditioned on the event $E$.

Let $\{\bar{z}_{k} \}$ and  $\{\bar{\mu}_{k}\}$ be the sequences of mean fields and strategies generated by Modified T-BR algorithm. We assume that Modified T-BR algorithm and TMFQ-learning algorithm have the same initialization, i.e.,  $\bar{z}_{0} = z_{0}$. Now, conditioned on the event $E$,
\begin{align*}
    \|\bar{z}_{k+1} - z_{k+1} \|_{1} & = \| \Phi(\bar{z}_{k}, \bar{\mu}_{k-1}) -  \Phi({z}_{k}, {\mu}_{k-1}) \|_{1} \\
    & \stackrel{(a)} = (1-\zeta) \sum_{s'}  |\sum_{s \in \mathcal{S}} \bar{z}_{k}(s) P_{\bar{z}_{k}, \bar{\mu}_{k-1}}(s'|s) - \sum_{s \in \mathcal{S}} {z}_{k}(s) P_{{z}_{k}, {\mu}_{k-1}}(s'|s) | + \zeta \sum_{s'} | \Psi(s') - \Psi(s') | \\
    & \stackrel{(b)} \leq (1-\zeta) (1+C_{2}) \|\bar{z}_{k} - {z}_{k} \|_{1} + (1-\zeta) C_{3} \|Q^{*}_{\bar{z}_{k-1}} - Q^{*}_{{k-1}} \|_{\infty} \\
    & \leq (1-\zeta) (1+C_{2}) \|\bar{z}_{k} - {z}_{k} \|_{1} + (1-\zeta) C_{3} \|Q^{*}_{\bar{z}_{k-1}} - Q^{*}_{z_{k-1}} \|_{\infty} + (1-\zeta) C_{3} \|Q^{*}_{z_{k-1}} - Q^{*}_{k-1} \|_{\infty} \\
    & \stackrel{(c)} \leq (1-\zeta) (1+C_{2}) \|\bar{z}_{k} - {z}_{k} \|_{1} + (1-\zeta) C_{3} D \|\bar{z}_{k-1} - z_{k-1} \|_{1} + (1-\zeta) C_{3} \epsilon_{3}
\end{align*}
Here $(a)$ follows from $\zeta$ regeneration event where $\Psi$ is the probability measure of the player regeneration process, $(b)$ follows from Lemma \ref{lem:mf-evol-by-Q} and $(c)$ follows from Lemma \ref{lem:Q-ub-z}.

Iteratively applying the above inequality, we get  $ \|\bar{z}_{k_{0}} - z_{k_{0}} \|_{1} \leq B \epsilon_{3},$ where $B = \frac{(2 (1-\zeta) A)^{k_{0}} (C_{3}\epsilon)}{(1-\zeta) (A-1)}$ where $A = \max \{ 1+C_{2}, C_{3}D \}$. Now, $\|\bar{z}^{*} - z_{k_{0}} \|_{1} \leq \|\bar{z}^{*} - \bar{z}_{k_{0}} \|_{1} + \|\bar{z}_{k_{0}} - z_{k_{0}} \|_{1} \leq \bar{\epsilon} + B \epsilon_{3}$ because $\|\bar{z}^{*} - \bar{z}_{k_{0}} \|_{1} \leq \bar{\epsilon}$ by the definition of $k_{0}$.

So, $\mathbb{P}( \|\bar{z}^{*} - z_{k_{0}} \|_{1} \leq \bar{\epsilon} + B \epsilon_{3}) \geq \mathbb{P}(E) = 1 - 2 k_{0} \delta$.

Setting $ \epsilon_{3} = \bar{\epsilon}/B$ and $\delta_{3} = \delta/2k_{0}$, and using the corresponding $I_{0}$ from \eqref{eq:synch-Q-sample-1}, we get the desired result.
\end{proof}

\section{Experiments}

\subsection{Parameters for Infection Spread Model}
We use the following parameters for simulations
 \begin{align*}
 \mathcal{|S|}&=25 \quad \mathcal{|A|}=5 \quad k=0.05 \quad \delta_1=1 \quad \delta_2=0.2 \quad &\delta_3=0.01 \quad  \zeta=0.1 \quad \epsilon=0.3 \quad \gamma =0.75   \\
 w_1 &\sim \begin{cases}
 \text{uniform} \{1,2,3\} \quad &\text{w.p} \quad 0.9 \\
 0 & \text{w.p} \quad 0.1
 \end{cases} 
 &w_2 \sim \begin{cases}
 \text{uniform} \{0, \cdots,s\} \quad &\text{w.p} \quad 0.9 \\
  0 & \text{w.p} \quad 0.1
 \end{cases}
 \end{align*}
For GMBL, we set $n_0=500$ and we run the outer loop for $500$ iterations. For TMFQL we run Q Learning for $1000$ time steps, and preform $5000$ iterations of the outer loop. For O-TMFQ Learning, run the outer loop for $5000$ iterations with different number of agents. For better sample efficiency, in both TMFQL and O-TMFQ Learning we initialize the Q function at each iteration with Q of the previous iteration. We used a logarithmically decaying learning rate, we decay the learning rate from $10^{-3}$ to $10^{-2}$.

For IQL we initialize each agent with a Q function and use the same parameters as O-TMFQ. We simulate a variant of MFQ where each agent estimates the meanfield to be the average state of a subset of the population and uses it to parameterize its Q function. We define this subset to be 512 agents chosen at random and kept constant for the duration of the simulation. Each agent also obtains samples from these 512 agents to updates its Q-Function. The other parameters are same as O-TMFQ.

\subsection{Amazon Mechanical Turk (MTurk) and other Gig Economy Marketplaces}
\begin{figure}
\centering
\begin{minipage}{.40\textwidth}
\centering
\includegraphics[width=1\columnwidth]{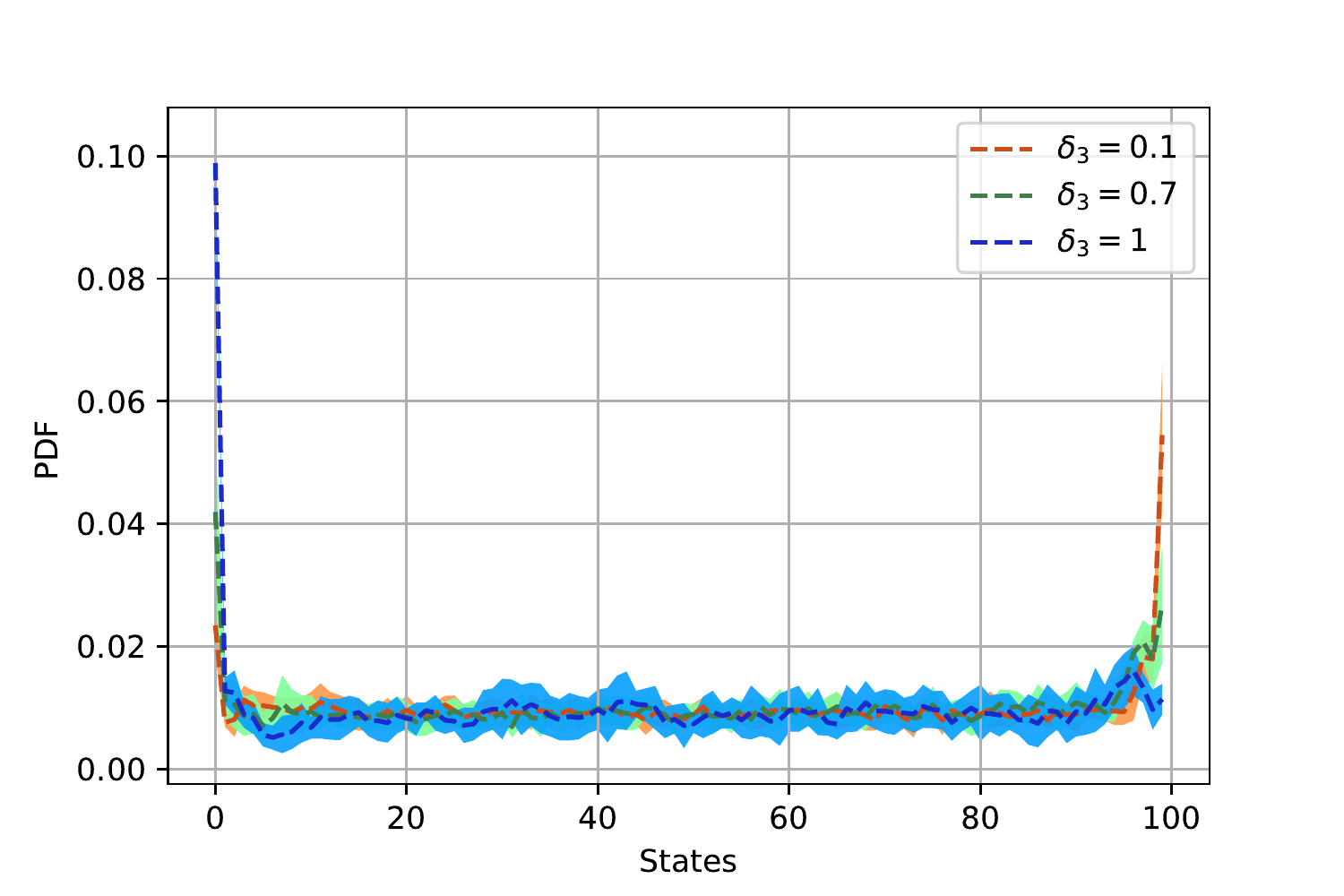}
\caption{PDF of O-TMFQ-Learning }
\label{fig:Turker_PDF_Algos}
\end{minipage}\hfill
\begin{minipage}{.40\textwidth}
\centering
\includegraphics[width=1\columnwidth]{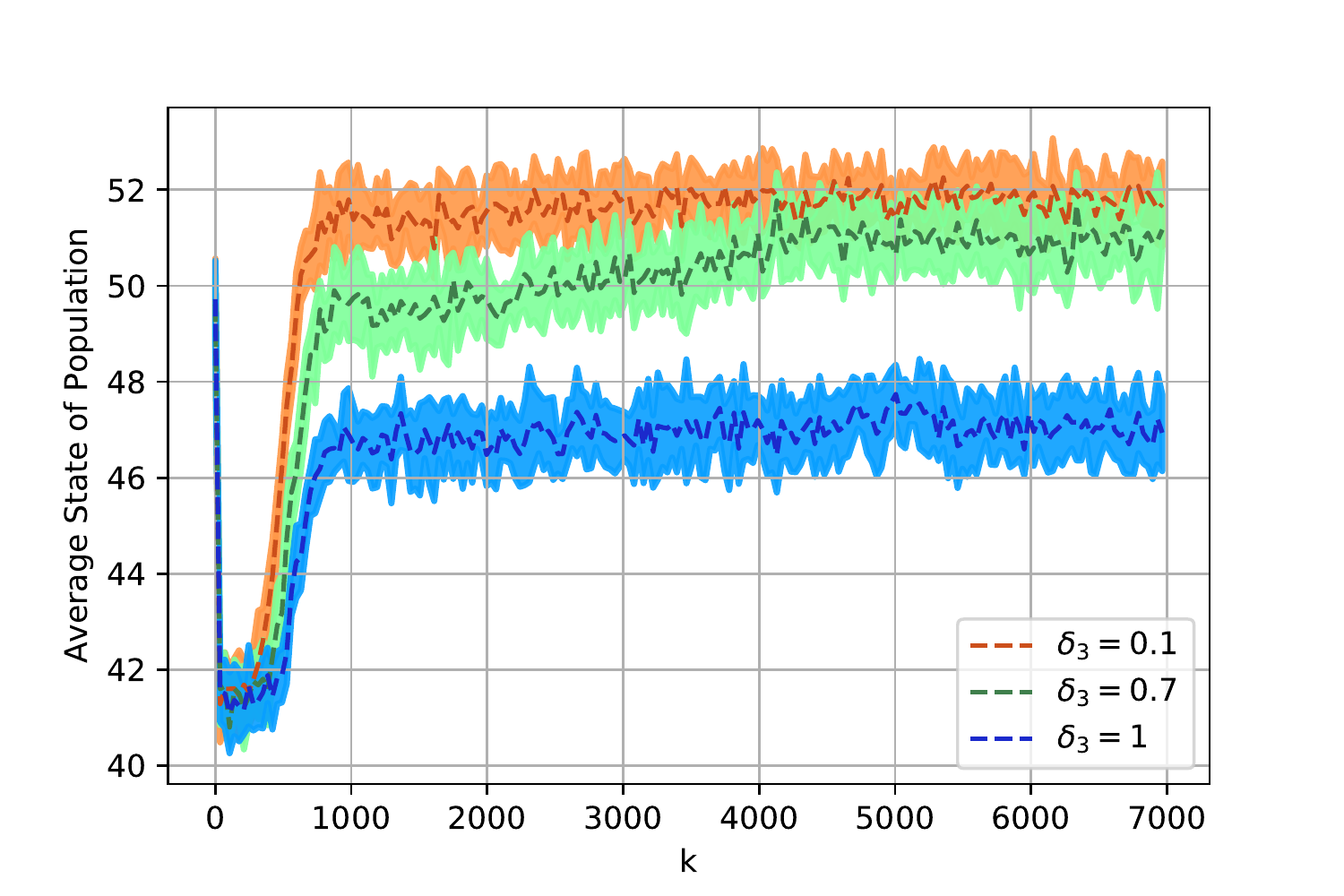}
\caption{Mean state of O-TMFQ-Learning}
\label{fig:Turker_mean_state_evol}
\end{minipage}\hfill
\begin{minipage}{.40\textwidth}
\centering
\includegraphics[width=1\columnwidth]{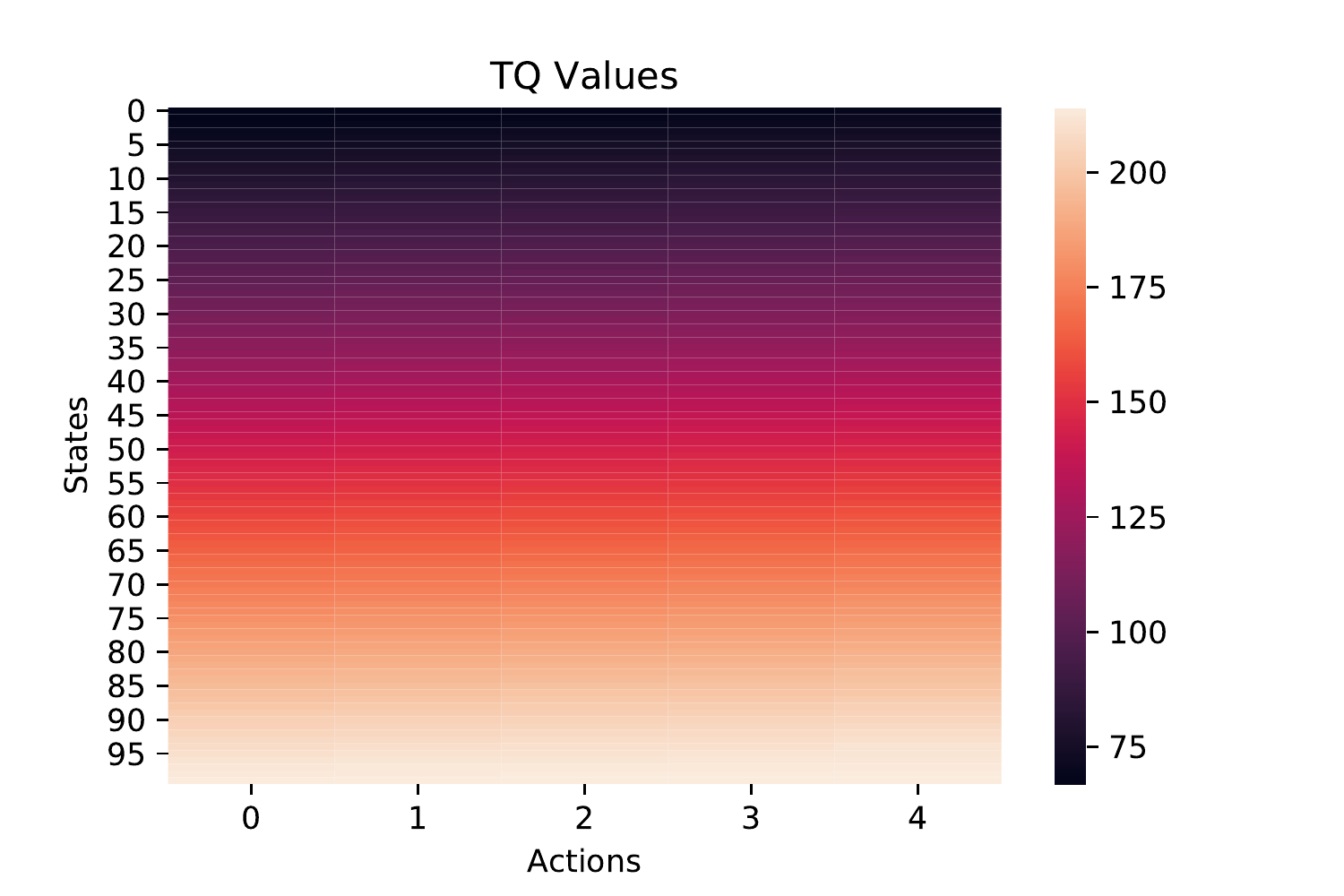}
\caption{Heat map of TQ Values for  $\delta_3=0.1.$ }
\label{fig:Turker_Q_heat}
\end{minipage}\hfill
\begin{minipage}{.40\textwidth}
\centering
\includegraphics[width=1\columnwidth]{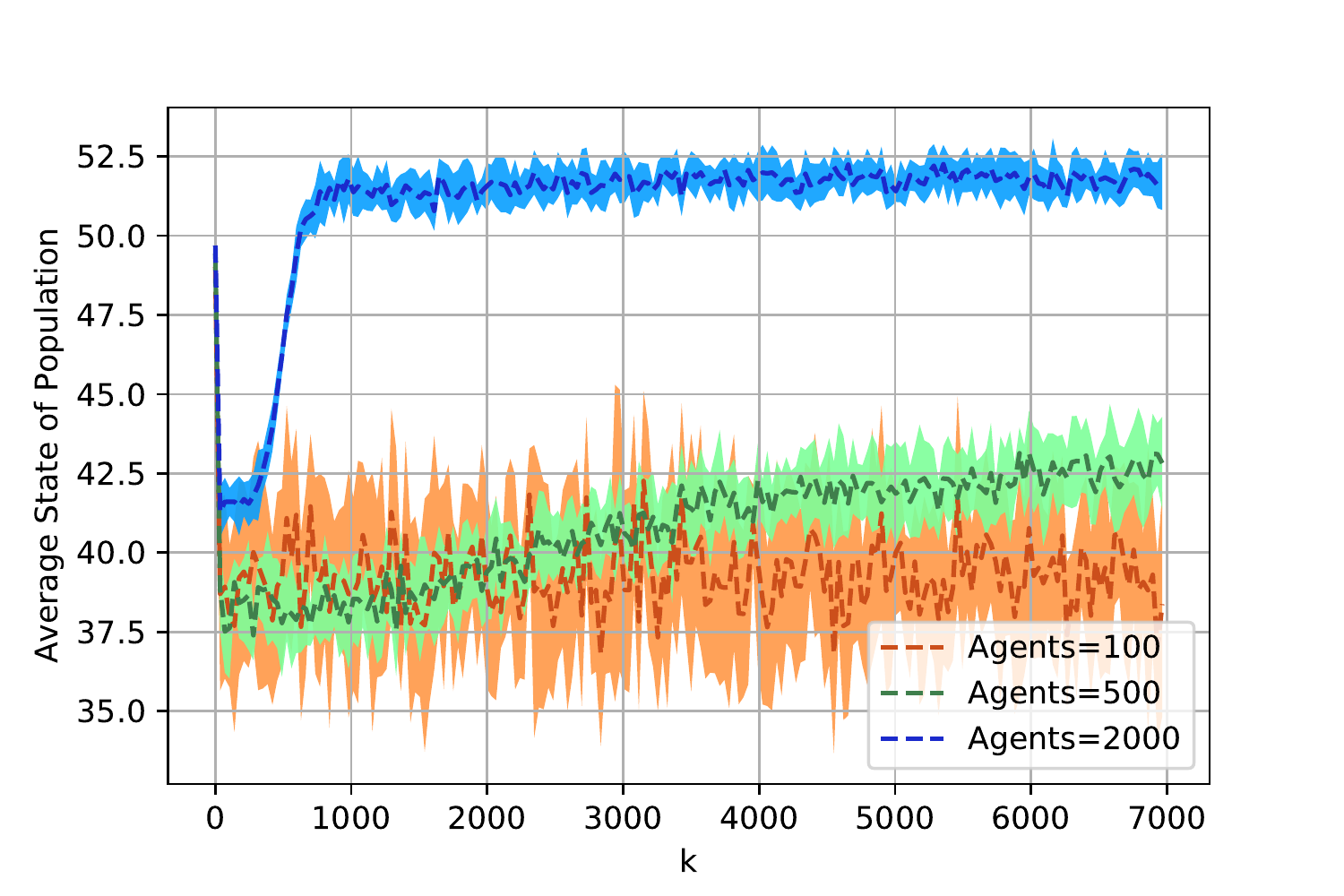}
\caption{Mean state of O-TMFQ-Learning: $\delta_3=0.1$}
\label{fig:Turker_mean_state_agents}
\end{minipage}\hfill
\end{figure}
We consider Amazon Mechanical Turk (MTurk) as an example of a Gig economy marketplace, which includes firms like Uber and Airbnb.   MTurk is a crowd sourcing market, wherein human workers are recruited to perform so-called Human Intelligence Tasks (HITs).  These HITs may take the form of labeling data sets or other tasks that are simple from a human's perspective, but might be difficult for machine learning to directly undertake.  The workers are called Turkers, and each has a quality score that depends on previous HITs undertaken.  The firm that originates these HITs may specify the price that it is willing to pay, as well as the minimum quality of the Turkers that it desires.  

There is a natural alignment of effort employed by Turkers in MTurk, since higher efforts translate into HITs done right, which in turn results in a higher quality, which finally results in firms willing to pay more per HIT.  Thus, if mean field quality is high, there is an incentive to perform HITs well and enhance ones' own quality.  This notion of incentive alignment applies to essentially all Gig economy marketplaces---the reputation of the agent directly enhances its reward, while the reputation of the marketplace as a whole (i.e., its mean field) draws customers willing to pay into the system, and so enhances the reward of the agent.

The formal system description is analogous to Infection Spread considered earlier, but we focus here on the application scenario. Thus, we have that each agent (Turker) has his/her quality state, and the strategic action is the choice of how much work to put into a HIT assigned to that Turker.  Higher effort implies higher cost, but also implies a higher improvement in the quality.  The overall reward is a combination that depends on the Turker's quality as well as the mean field quality. 

Let $s$ denote the quality of a Turker's profile, and let action $a$ denote the effort the turker puts in to maintain the quality of the profile; this may include number of jobs successfully completed, time taken to complete a job etc. Let $c(a)=\delta_3a$ denote the cost incurred in performing action $a$. The quality perceived by an entity offering jobs depends on both the quality of the individual Turker and the population as a whole (via the mean field). Thus, the reward to a Turker is a function of the perceived quality and cost incurred in taking action. We define state transition and reward as follows, 
\begin{align*}
    s'& = (s+a-w_1)_{+}\mathbbm{1}\{E_1\} + w_2 \mathbbm{1}\{E_2\} \\
    r &= \delta_1s+\delta_2 \sum_{s\in \mathcal{S}}sz(s) -\delta_3 a,
\end{align*}
where $E_1,E_2$ are mutually exclusive events that occur with probabilities $1-\zeta, \zeta,$ respectively, and $w_1,w_2$ are realizations of non-negative integer random variables.

We use the following parameters for simulations, 
\begin{align*}
|\mathcal{S}|&=100 \quad |\mathcal{A}|=5 \quad \delta_1=0.5 \quad \delta_2=0.2 \quad \zeta = 0.1 \quad \epsilon=0.3 \quad \gamma=0.75 \\
w_1 &\sim \text{uniform} \{0,1,2,3\}   \quad w_2 \sim \begin{cases}
 \text{uniform} \{0, \cdots,|\mathcal{S}|\} \quad &\text{w.p} \quad 0.9 \\
  0 & \text{w.p} \quad 0.1
 \end{cases}
\end{align*}
For O-TMFQ Learning, we run the outer iteration for $7000$ steps. As before, for better sample efficiency we initialize the Q function at each iteration with Q of the previous iteration. We used a logarithmically decaying learning rate, we decay the learning rate from $10^{-3}$ to $10^{-2}$.

The behavior of our RL algorithms is much the same as the earlier case, and is shown in Figures~\ref{fig:Turker_PDF_Algos}--\ref{fig:Turker_mean_state_agents}. Figure \ref{fig:Turker_PDF_Algos} shows the pdf of the final mean field distribution obtained by performing O-TMFQL with $2000$ agents with different values of $\delta_3$ while figure \ref{fig:Turker_mean_state_evol} shows the average state of the population. Observe that as the cost of action $\delta_3$ increases, agents take lower actions and are hence distributed towards lower states. Figure \ref{fig:Turker_Q_heat} is a heat map of the final TQ value function. Figure \ref{fig:Turker_mean_state_agents} shows the mean state evolution with different number of agents for  $\delta_3=0.1$. Observe that the convergence is poor with lesser number of agents.

\end{document}